\documentclass[twoside,a4paper]{amsart}
\providecommand{\text}[1]{\mbox{#1}}
\usepackage[OT2, OT1]{fontenc}
\usepackage{animate}
\usepackage{listings}
 \usepackage[pdfauthor={Vladimir V. Kisil},%
      pdftitle={Lectures on Moebius-Lie geometry and its extension},%
      pdfsubject={mathematics, physics},%
      backref=page,%
    pdfkeywords={symmetries, quantum mechanics, hypercomplex numbers},
  breaklinks=true,%
  colorlinks=true,%
  bookmarks=true,%
  pagebackref=true
  ]{hyperref}
\usepackage[nobysame, initials,bibtex-style]{amsrefs}

\providecommand{\citelist}[1]{#1}

\makeatletter
  \providecommand{\limlike}[1]{\mathop {\operator@font #1}}
  \providecommand{\loglike}[1]{\mathop {\operator@font #1}\nolimits}
\makeatother

\makeatletter
\renewcommand{\PrintNames@a}[4]{%
    \PrintSeries{\name}
        {#1}
        {}{ and \set@othername}
        {,}{ \set@othername}
        {}{ and \set@othername}
        {#2}{#4}{#3}%
}
\makeatother

\BibSpecAlias{incollection}{inproceedings}
\BibSpec{article}{%
    +{}  {\PrintAuthors}                {author}
    +{,} { \emph}                            {title}
    +{,} { }                            {part}
    +{:} { }                            {subtitle}
    +{,} { \PrintContributions}         {contribution}
    +{,} { \PrintPartials}              {partial}
    +{,} { }                       {journal}
    +{} { \textbf}                            {volume}
    +{} { \parenthesize}                 {date}
    +{} { }                             {pages}
    +{,} { }                            {status}
    +{,} { \PrintTranslation}           {translation}
    +{,} { Reprinted in \PrintReprint}  {reprint}
    +{,} { \doi}  {doi}
    +{,} { \url}  {url}
    +{,} { }                            {note}
    +{.} {}                             {transition}
}

\BibSpec{partial}{%
    +{}  {}                             {part}
    +{:} { }                            {subtitle}
    +{,} { \PrintContributions}         {contribution}
    +{,} { \emph}                       {journal}
    +{,} { \textbf}                            {volume}
    +{}  { \parenthesize}               {number}
    +{:} {}                             {pages}
    +{.} { \PrintDateB}                 {date}
}

\BibSpec{book}{%
    +{}  {\PrintPrimary}                {transition}
    +{,} { \emph}                       {title}
    +{,} { }                            {part}
    +{:} { \emph}                       {subtitle}
    +{,} { }                            {series}
    +{,} { \voltext}                    {volume}
    +{,} { Edited by \PrintNameList}    {editor}
    +{,} { Translated by \PrintNameList}{translator}
    +{,} { \PrintContributions}         {contribution}
    +{,} { \PrintEdition}               {edition}
    +{,} { }                            {publisher}
    +{,} { }                            {address}
    +{ } { \PrintDateB}                 {date}
    +{.} {}                             {transition}
}

\BibSpec{collection.article}{%
    +{}  {\PrintAuthors}                {author}
    +{,} { \emph}                            {title}
    +{,} { }                            {part}
    +{:} { }                            {subtitle}
    +{,} { \PrintContributions}         {contribution}
    +{,} { \PrintConference}            {conference}
    +{,} { \PrintBook}                  {book}
    +{,} { In \PrintEditorsA}              {editor}
    +{,} { }                         {booktitle}
    +{,} { pages~}                      {pages}
    +{,} { }                            {publisher}
    +{,} { }                            {organization}
    +{,} { }                            {address}
    +{,} { \PrintDateB}                 {date}
    +{,} { \PrintTranslation}           {translation}
    +{,} { Reprinted in \PrintReprint}  {reprint}
    +{,} { }                            {note}
    +{.} {}                             {transition}
}

 \BibSpec{inproceedings}{%
     +{}  {\PrintAuthors}                {author}
     +{,} { \emph}                            {title}
     +{,} { }                            {part}
     +{:} { }                            {subtitle}
     +{,} { \PrintContributions}         {contribution}
     +{,} { \PrintConference}            {conference}
     +{,} { \PrintBook}                  {book}
     +{,} { In \PrintEditorsA}              {editor}
     +{,} {  }                         {booktitle}
     +{,} { pages~}                      {pages}
     +{,} { }                            {publisher}
     +{,} { }                            {organization}
     +{,} { }                            {address}
     +{,} { \PrintDateB}                 {date}
     +{,} { \PrintTranslation}           {translation}
     +{,} { Reprinted in \PrintReprint}  {reprint}
     +{,} { }                            {note}
     +{.} {}                             {transition}
 }

\BibSpec{conference}{%
    +{}  {}                        {title}
    +{}  {\PrintConferenceDetails} {transition}
}

\BibSpec{innerbook}{%
    +{,} { \emph}                       {title}
    +{,} { }                            {part}
    +{:} { \emph}                       {subtitle}
    +{,} { }                            {series}
    +{,} { \voltext}                    {volume}
    +{,} { Edited by \PrintNameList}    {editor}
    +{,} { Translated by \PrintNameList}{translator}
    +{,} { \PrintContributions}         {contribution}
    +{,} { }                            {publisher}
    +{,} { }                            {organization}
    +{,} { }                            {address}
    +{,} { \PrintEdition}               {edition}
    +{,} { \PrintDateB}                 {date}
    +{,} { }                            {note}
    +{.} {}                             {transition}
}

\BibSpec{report}{%
    +{}  {\PrintPrimary}                {transition}
    +{,} { \emph}                       {title}
    +{,} { }                            {part}
    +{:} { \emph}                       {subtitle}
    +{,} { \PrintContributions}         {contribution}
    +{,} { Technical Report }           {number}
    +{,} { }                            {series}
    +{,} { }                            {organization}
    +{,} { }                            {address}
    +{,} { \PrintDateB}                 {date}
    +{,} { \PrintTranslation}           {translation}
    +{,} { Reprinted in \PrintReprint}  {reprint}
    +{,} { }                            {note}
    +{.} {}                             {transition}
}

\BibSpec{thesis}{%
    +{}  {\PrintAuthors}                {author}
    +{,} { \emph}                       {title}
    +{:} { \emph}                       {subtitle}
    +{,} { \PrintThesisType}            {type}
    +{,} { }                            {organization}
    +{,} { }                            {address}
    +{,} { \PrintDateB}                 {date}
    +{,} { \PrintTranslation}           {translation}
    +{,} { Reprinted in \PrintReprint}  {reprint}
    +{,} { }                            {note}
    +{.} {}                             {transition}
}
\usepackage{graphicx,amssymb}

\makeatletter
  \providecommand{\limlike}[1]{\mathop {\operator@font #1}}
  \providecommand{\loglike}[1]{\mathop {\operator@font #1}\nolimits}
\makeatother

\providecommand{\nscalar}[3][\relax]{\left[ #2,#3
        \right]\ifx#1\relax\else_{#1}\fi}
\providecommand{\tr}{\mathop{\mathrm{tr}}}
\newcommand{\CPP}{\textsf{C++}}

\newcommand{\Python}{\textsf{Python}}
\newcommand{\NoWEB}{\texttt{noweb}}
\providecommand{\MetaPost}{\texttt{Meta}\-\texttt{Post}}
\providecommand{\GiNaC}{\textsf{GiNaC}}
\providecommand{\pyGiNaC}{\textsf{pyGiNaC}}
\providecommand{\href}[2]{#2}
\providecommand{\SL}[1][2]{\ensuremath{\mathrm{SL}_{#1}(\Space{R}{})}}

\newcommand{\Aprime}{A'}
\providecommand{\Space}[3][]{\ensuremath{\mathbb{#2}^{#3}_{#1}{}}}
  \providecommand{\FSpace}[3][]{\ensuremath{\ifx#2l \ell_{#3}^{#1}{}\else
  \mathcal{#2}_{#3}^{#1}{}\fi}} 
\providecommand{\norm}[2][\relax]{\left\|#2\right\|\ifx#1\relax\else_{#1}\fi}
\providecommand{\modulus}[2][\relax]{\left| #2 \right|\ifx#1\relax\else_{#1}\fi}
\providecommand{\scalar}[3][\relax]{\left\langle #2,#3 
        \right\rangle\ifx#1\relax\else_{#1}\fi}

\providecommand{\alli}{\iota}

\providecommand{\notingiq}{}
\providecommand{\notingiqsemi}{}
\providecommand{\notingiqsemitocoma}{,}
\providecommand{\rmi}{\mathrm{i}}

\providecommand{\map}[1]{\mathsf{#1}}
\providecommand{\spec}[1][]{\ensuremath{\mathbf{sp}}\,}

\providecommand{\cycle}[3][\relax]{\ifx#1\relax C\else{\ifx#1\tilde S
    \else G \fi}\fi^{#2}_{#3}}
\providecommand{\Cliff}[2][\comment]{{\ensuremath{%
\mathcal{C}\kern-0.18em\ell(#1,#2)}}}
\providecommand{\comment}[1]{}
\providecommand{\uir}[3][0]{\ifcase #1{\rho^{#2}_{#3}}%
\or {\breve{\rho}^{#2}_{#3}}%
\or {\tilde{\rho}^{#2}_{#3}}\fi}

  \providecommand{\Zbl}[1]{Zbl\href{http://www.emis.de:80/cgi-bin/zmen/ZMATH/en/zmathf.html?first=1&maxdocs=3&type=html&an=#1&format=complete}{#1}}
\providecommand{\bsigma}{\tau}

\providecommand{\myeprint}[2]{E-print: \href{#1}{\texttt{#2}}}
\providecommand{\doi}[1]{doi: \href{http://dx.doi.org/#1}{#1}}

\providecommand{\oper}[1]{\mathcal{#1}}
\newcommand*\vtick{\kern -.1em\textsc{\char13}}

\iftrue
\newtheorem{theorem}{Theorem}[section]

    \PassOptionsToPackage{british}{babel}
    \newtheorem{proposition}[theorem]{Proposition}
    \newtheorem{lemma}[theorem]{Lemma}
    \newtheorem{corollary}[theorem]{Corollary}
    
    \theoremstyle{definition}

    \newtheorem{definition}[theorem]{Definition}

    \newtheorem{example}[theorem]{Example}

    \theoremstyle{remark}
    \newtheorem{remark}[theorem]{Remark}
   \providecommand{\such}{\,\mid\,}
   \DeclareMathOperator{\arccosh}{arccosh}
\else
   \providecommand{\such}{\,\semicolon\,}
\fi

\theoremstyle{remark}
\newtheorem{claim}[theorem]{Claim}
\providecommand{\dSpace}[2]{\dot{\mathbb{#1}}^{#2}}

\makeindex

\usepackage[all]{xy}
\DeclareFontFamily{U}{mathx}{\hyphenchar\font45}
\DeclareFontShape{U}{mathx}{m}{n}{
      <5> <6> <7> <8> <9> <10>
      <10.95> <12> <14.4> <17.28> <20.74> <24.88>
      mathx10
      }{}
\DeclareSymbolFont{mathx}{U}{mathx}{m}{n}
\DeclareFontSubstitution{U}{mathx}{m}{n}
\DeclareMathAccent{\wideparen}{0}{mathx}{"75}
\DeclareFontFamily{OT1}{cyr}{}
\DeclareFontShape{OT1}{cyr}{m}{n}
   {  <5> <6> <7> <8> <9> gen * wncyr
      <10> <10.95> <12> <14.4> <17.28> <20.74> <24.88> wncyr10}{}
\DeclareFontShape{OT1}{cyr}{m}{it}
    {
       <5> <6> <7> <8> <9> gen * wncyi
      <10> <10.95> <12> <14.4> <17.28> <20.74> <24.88>wncyi10
      }{}
\DeclareFontShape{OT1}{cyr}{m}{ss}
    {
       <5> <6> <7> <8> wncyss8
       <9> wncy9
      <10> <10.95> <12> <14.4> <17.28> <20.74> <24.88>wncyss10
      }{}
\DeclareFontShape{OT1}{cyr}{m}{sc}
    {
       <5> <6> <7> <8> <9> <10> <10.95> <12> <14.4> <17.28> <20.74> <24.88>wncysc10
      }{}
\DeclareFontShape{OT1}{cyr}{bx}{n}
   {
       <5> <6> <7> <8> <9> gen * wncyb
      <10> <10.95> <12> <14.4> <17.28> <20.74> <24.88>wncyb10
      }{}
\DeclareTextFontCommand{\textcyr}{\fontfamily{cyr}\selectfont}
\newcommand{\cyr}{\fontfamily{cyr}\selectfont\def\cprime{\~}}
\providecommand{\cprime}{'}

\begin{document}

\title[M\"obius--Lie Geometry and Its Extension]
{Lectures on M\"obius--Lie Geometry\\
  and Its Extension}

\author[Vladimir V. Kisil]%
{Vladimir V. Kisil}

\address{School of Mathematics, University of Leeds,  Leeds, LS2\,9JT, England}
\email{email: \href{mailto:kisilv@maths.leeds.ac.uk}{kisilv@maths.leeds.ac.uk}\\}
\urladdr{URL: \href{http://www.maths.leeds.ac.uk/~kisilv/}%
{http://www.maths.leeds.ac.uk/\~{}kisilv/}}


\begin{abstract}
  These lectures review the classical M\"obius--Lie geometry and
  recent work on its extension. The latter considers ensembles of
  cycles (quadrics), which are interconnected through
  conformal-invariant geometric relations (e.g. ``to be orthogonal'',
  ``to be tangent'', etc.), as new objects in an extended
  M\"obius--Lie geometry.  It is shown on examples, that such
  ensembles of cycles naturally parameterise many other
  conformally-invariant families of objects, two examples---the
  Poincar\'e extension and continued fractions are considered in
  detail. Further examples, e.g. loxodromes, wave fronts and
  integrable systems, are published elsewhere.

  The extended M\"obius--Lie geometry is efficient due to a method,
  which reduces a collection of conformally in\-vari\-ant geometric
  relations to a system of linear equations, which may be accompanied
  by one fixed quadratic relation. The algorithmic nature of the
  method allows to implement it as a {\CPP} library, which operates
  with numeric and symbolic data of cycles in spaces of arbitrary
  dimensionality and metrics with any signatures. Numeric calculations
  can be done in exact or approximate arithmetic. In the two- and
  three-dimensional cases illustrations and animations can be
  produced. An interactive {\Python} wrapper of the library is
  provided as well.
  \end{abstract}
\keywords{M\"obius--Lie geometry, spheres geometry, Poincar\'e
  extension, continued fraction, integrable system, loxodrome,
  fraction-linear transformation, Clifford algebra, indefinite inner
  product space.}
\subjclass[2000]{Primary 51B25; Secondary 30B70, 51M05,  51N25, 51B10, 68U05, 11E88, 68W30.}

\maketitle

\tableofcontents

\section{Introduction}

\href{https://en.wikipedia.org/wiki/Lie_sphere_geometry}{Lie sphere
  geometry}~\citelist{ \cite{Cecil08a} \cite{Benz07a}*{Ch.~3}} in the
simplest planar setup unifies circles, lines and points---all together
called \emph{cycles} in this setup.  Symmetries of Lie spheres
geometry include (but are not limited to) fractional linear
transformations (FLT) of the form:
\begin{equation}
  \label{eq:flt-defn}
  \begin{pmatrix}
    a&b\\c&d
  \end{pmatrix}:\  x \mapsto
  \frac{ax+b}{cx+d}\,, \qquad \text{where }
  \det\begin{pmatrix}
    a&b\\c&d
  \end{pmatrix}\neq 0.
\end{equation}
Following other sources, e.g. ~\cite{Simon11a}*{\S~9.2}, we call
\eqref{eq:flt-defn} by FLT and reserve the name ``M\"obius maps'' for
the subgroup of FLT which fixes a particular cycle. For example, on
the complex plane FLT are generated by elements of
\(\mathrm{SL}_2(\Space{C}{})\) and M\"obius maps fixing the real line
are produced by
\(\mathrm{SL}_2(\Space{R}{})\)~\cite{Kisil12a}*{Ch.~1}.

There is a natural set of FLT-invariant geometric relations between
cycles (to be orthogonal, to be tangent, etc.) and the restriction of
Lie sphere geometry to invariants of FLT is called \emph{M\"obius--Lie
  geometry}.  Thus, an ensemble of cycles, structured by a set of such
relations, will be mapped by FLT to another ensemble with the same
structure.

It was shown recently that ensembles of cycles with certain
FLT-invariant relations provide helpful parametrisations of new
objects, e.g. points of the Poincar\'e extended space~\cite{Kisil15a},
loxodromes~\cite{KisilReid18a} or continued
fractions~\cites{BeardonShort14a,Kisil14a}, see
Example~\ref{ex:ensamble-math} below for further details. Thus, we
propose \emph{to extend M\"obius--Lie geometry and consider ensembles
  of cycles as its new objects}, cf. formal
Defn.~\ref{de:extended-Lie-Mobius}. Naturally, ``old''
objects---cycles---are represented by simplest one-element ensembles
without any relation. This paper provides conceptual foundations of
such extension and demonstrates its practical implementation as a
{\CPP} library {{\bf{}figure}}\footnote{All described software is licensed
  under GNU GPLv3~\cite{GNUGPL}.}. Interestingly, the development of
this library shaped the general approach, which leads to specific
realisations in~\cites{Kisil15a,Kisil14a,KisilReid18a}.

More specifically, the library {{\bf{}figure}} manipulates ensembles of
cycles (quadrics) interrelated by certain FLT-invariant geometric
conditions.  The code is build on top of the previous library
{{\bf{}cycle}}~\cites{Kisil05b,Kisil12a,Kisil06a}, which manipulates
individual cycles within the \GiNaC~\cite{GiNaC} computer algebra
system. Thinking an ensemble as a graph, one can say that the library
{{\bf{}cycle}} deals with individual vertices (cycles), while {{\bf{}figure}}
considers edges (relations between pairs of cycles) and the whole
graph. Intuitively, an interaction with the library {{\bf{}figure}} reminds
compass-and-straightedge constructions, where new lines or circles are
added to a drawing one-by-one through relations to already presented
objects (the line through two points, the intersection point or the
circle with given centre and a point). See
Example~\ref{ex:touch-centres-collinear} of such interactive
construction from the {\Python} wrapper, which provides an analytic
proof of a simple geometric statement.

It is important that both libraries are capable to work in spaces of
any dimensionality and metrics with an arbitrary signatures:
Euclidean, Minkowski and even degenerate. Parameters of objects can be
symbolic or numeric, the latter admit calculations with exact or
approximate arithmetic.  Drawing routines work with any (elliptic,
parabolic or hyperbolic) metric in two dimensions and the euclidean
metric in three dimensions.

The mathematical formalism employed in the library {{\bf{}cycle}} is based
on Clifford algebras, which are intimately connected to fundamental
geometrical and physical objects
\cites{HestenesSobczyk84a,Hestenes15a}. Thus, it is not surprising
that Clifford algebras have been already used in various geometric
algorithms for a long time, for example see
\cites{Hildenbrand13a,Vince08a,DorstDoranLasenby02a} and further
references there. Our package deals with cycles through
Fillmore--Springer--Cnops construction (FSCc) which also has a long
history, see~\citelist{\cite{Schwerdtfeger79a}*{\S~1.1}
  \cite{Cnops02a}*{\S~4.1}
  \cite{FillmoreSpringer90a} \cite{Kirillov06}*{\S~4.2}
  \cite{Kisil05a} \cite{Kisil12a}*{\S~4.2}} and
section~\ref{sec:lie-spheres-geometry} below. Compared to a plain
analytical treatment~\citelist{\cite{Pedoe95a}*{Ch.~2}
  \cite{Benz07a}*{Ch.~3}}, FSCc is much more efficient and
conceptually coherent in dealing with FLT-invariant properties of
cycles. Correspondingly, the computer code based on FSCc is easy to
write and maintain.

The paper outline is as follows. In Section~\ref{sec:math-backgr} we
sketch the mathematical theory (M\"obius--Lie geometry) covered by the
package of the previous library {{\bf{}cycle}}~\cite{Kisil05b} and the
present library {{\bf{}figure}}. We expose the subject with some
references to its history since this can facilitate further
development. Sec.~\ref{sec:conn-quadr-cycl} describes the principal
mathematical tool used by the library {{\bf{}figure}}.  It allows to
reduce a collection of various linear and quadratic equations
(expressing geometrical relations like orthogonality and tangency) to
a set of linear equations and \emph{at most one} quadratic
relation~\eqref{eq:det-normalisation-cond}. Notably, the quadratic
relation is the same in all cases, which greatly simplifies its
handling. This approach is the cornerstone of the library
effectiveness both in symbolic and numerical computations.

We consider two examples of ensembles of cycles in
details. Section~\ref{sec:poinc-extens-mobi} presents the Poincar\'e
extension. Given sphere preserving (M\"obius) transformations in
\(n\)-dimensional Eucli\-delta\-an space one can use the Poincar\'e
extension to obtain sphere preserving transformations in a half-space
of \(n+1\) dimensions. The Poincar\'e extension is usually provided
either by an explicit formula or by some geometric construction. Here
we present its algebraic background and describe all available
options. The solution is given either in terms of one-parameter
subgroups of M\"obius transformations or ensembles of cycle
representing equivalent triples of quadratic forms.
  
Another example, presented in Section~\ref{sec:exampl-cont-fract},
uses interrelations between continued fractions, M\"obius
transformations and representations of cycles by \(2\times 2\)
matrices. This leads us to several descriptions of continued fractions
through ensembles of cycles consisting of chains of orthogonal or
touching horocycles. One of these descriptions was proposed in a
recent paper by A.~Beardon and I.~Short. The approach is extended to
several dimensions in a way which is compatible to the early
propositions of A.~Beardon based on Clifford algebras.

In Sec.~\ref{sec:figures-as-families} we present several examples of
ensembles, which were already used in mathematical
theories~\cites{Kisil15a,Kisil14a,KisilReid18a}, then we describe how
ensembles are encoded in the present library {{\bf{}figure}} through
the functional programming framework.

Sec.~\ref{sec:mathematical-results} outlines several typical usages of
the package. An example of a new statement discovered and demonstrated
by the package is given in Thm.~\ref{th:nine-points}. 

All coding-related material is enclosed as appendices in the full
documentation on the project page~\cite{Kisil05b}. They contain:
\begin{enumerate}
\item Numerous examples of the library usage
  starting from the very simple ones.
\item  A systematic list of callable
  methods.
\item Actual code of the library.
\end{enumerate}
Sec.~\ref{sec:math-backgr}, Example~\ref{ex:touch-centres-collinear}
below or the above-mentioned first two appendices of the full
documentation can serve as an entry point for a reader with respective
preferences and background.

\section[Moebius--Lie Geometry and the cycle Library]{M\"obius--Lie Geometry and the {{\bf{}cycle}} Library}
\label{sec:math-backgr}
We briefly outline mathematical formalism of the extend M\"obius--Lie
geometry, which is implemented in the present package. We do not aim
to present the complete theory here, instead we provide a minimal
description with a sufficient amount of references to further
sources. The hierarchical structure of the theory naturally splits the
package into two components: the routines handling individual cycles
(the library {{\bf{}cycle}} briefly reviewed in this section), which were
already introduced elsewhere~\cite{Kisil05b}, and the new component
implemented in this work, which handles families of interrelated
cycles (the library {{\bf{}figure}} introduced in the next section).

\subsection[Moebius--Lie geometry and FSC construction]{M\"obius--Lie geometry and FSC construction}
\label{sec:lie-spheres-geometry}
M\"obius--Lie geometry in \(\Space{R}{n}\) starts from an observation
that points can be treated as spheres of zero radius and planes are
the limiting case of spheres with radii diverging to
infinity. Oriented spheres, planes and points are called together
\emph{cycles}\index{cycle}. Then, the second crucial step is to treat
cycles not as subsets of \(\Space{R}{n}\) but rather as points of some
projective space of higher dimensionality,
see~\citelist{\cite{Benz08a}*{Ch.~3} \cite{Cecil08a} \cite{Pedoe95a}
  \cite{Schwerdtfeger79a}}.

To distinguish two spaces we will call \(\Space{R}{n}\) as the \emph{point
  space}%
\index{point!space}%
\index{space!point} and the higher dimension space, where cycles are
represented by points---the \emph{cycle space}%
\index{cycle!space}%
\index{space!cycle}. Next important observation is that geometrical
relations between cycles as subsets of the point space can be
expressed in term of some indefinite metric on the cycle
space. Therefore, if an indefinite metric shall be considered anyway,
there is no reason to be limited to spheres in Euclidean space \(\Space{R}{n}\)
only. The same approach shall be adopted for quadrics in
spaces \(\Space{R}{pqr}\) of an arbitrary signature
\(p+q+r=n\), including \(r\) nilpotent elements,
cf.~\eqref{eq:clifford-defn} below.

A useful addition to M\"obius--Lie geometry is provided by the
Fillmore--Springer--Cnops construction
(FSCc)~\citelist{\cite{Schwerdtfeger79a}*{\S~1.1}
  \cite{Cnops02a}*{\S~4.1} \cite{Porteous95}*{\S~18}
  \cite{FillmoreSpringer90a} \cite{Kirillov06}*{\S~4.2}
  \cite{Kisil05a} \cite{Kisil12a}*{\S~4.2}}. It is a correspondence
between the cycles (as points of the cycle space) and certain
\(2\times 2\)-matrices defined in~\eqref{eq:spheres-Rn} below. The
main advantages of FSCc are:
\begin{enumerate}
\item The correspondence between cycles and matrices respects the
  projective structure of the cycle space.
\item The correspondence is FLT covariant.
\item The indefinite metric on the cycle space can be expressed
  through natural operations on the respective matrices.
\end{enumerate}
The last observation is that for restricted groups of M\"obius
transformations the metric of the cycle space may not be completely
determined by the metric of the point space,
see~\citelist{\cite{Kisil06a} \cite{Kisil05a}
  \cite{Kisil12a}*{\S~4.2}} for an example in two-dimensional space.

FSCc is useful in consideration of the Poincar\'e extension of
M\"obius maps~\cite{Kisil15a}, loxodromes~\cite{KisilReid18a} and
continued fractions~\cites{Kisil14a}. In theoretical physics FSCc
nicely describes conformal compactifications of various space-time
models~\citelist{\cite{HerranzSantander02b} \cite{Kisil06b}
  \cite{Kisil12a}*{\S~8.1}}.  Regretfully, FSCc have not yet
propagated back to the most fundamental case of complex numbers,
cf.~\cite{Simon11a}*{\S~9.2} or somewhat cumbersome techniques used
in~\cite{Benz07a}*{Ch.~3}. Interestingly, even the founding fathers
were not always strict followers of their own techniques,
see~\cite{FillmoreSpringer00a}.

We turn now to the explicit definitions.

\subsection{Clifford algebras, FLT transformations, and Cycles}
\label{sec:cliff-algebr-mobi}
We describe here the mathematics behind the the first library called
{{\bf{}cycle}}, which implements fundamental geometrical relations
between quadrics in the space \(\Space{R}{pqr}\) with the
dimensionality \(n=p+q+r\) and metric
\(x_1^2+\ldots+x_p^2-x_{p+1}^2-\ldots-x_{p+q}^2\). A version
simplified for complex numbers only can be found
in~\cites{Kisil15a,KisilReid18a,Kisil14a}. In two dimensions usage of
dual and double numbers is preferable due to their
commutativity~\cites{Kisil17a,Kisil12a,Kisil06a,Yaglom79}.

The Clifford algebra \(\Cliff{p,q,r}\) is the associative unital algebra over
\(\Space{R}{}\) generated by the elements \(e_1\),\ldots,\(e_n\)
satisfying the following relation:
\begin{equation}
\label{eq:clifford-defn}
  e_i e_j =- e_je_i\,, \quad \text{ and } \quad e_i^2=\left\{
    \begin{array}{ll}
      -1,&\text{ if } 1\leq i\leq p{\notingiqsemi}\\
      1,&\text{ if } p+1\leq i\leq p+q{\notingiqsemi}\\
      0,&\text{ if } p+q+1\leq i\leq p+q+r.
    \end{array}
  \right.
\end{equation}
It is common~\cites{DelSomSou92,Cnops02a,Porteous95,%
  HestenesSobczyk84a,Hestenes15a} to consider mainly Clifford algebras
\(\Cliff{n}=\Cliff{n,0,0}\) of the Euclidean space or the algebra
\(\Cliff{p,q}=\Cliff{p,q,0}\) of the pseudo-Euclidean (Minkowski)
spaces. However, Clifford algebras \(\Cliff{p,q,r}\), \(r>0\) with
nilpotent generators \(e_i^2=0\) correspond to interesting
geometry~\cites{Kisil12a,Kisil05a,Yaglom79,Mustafa17a} and
physics~\cites{GromovKuratov06a,Gromov10a,Gromov12a,%
  Kisil12c,Kisil09e,Kisil17a} as well. Yet, the geometry with
idempotent units in spaces with dimensionality \(n>2\) is still not
sufficiently elaborated.

An element of
\(\Cliff{p,q,r}\) having the form \(x=x_1e_1+\ldots+x_ne_n\) can be
associated with the vector \((x_1,\ldots,x_n)\in\Space{R}{pqr}\).  The
\emph{reversion} \(a\mapsto a^*\) in
\(\Cliff{p,q,r}\)~\cite{Cnops02a}*{(1.19(ii))} is defined on vectors by
\(x^*=x\) and extended to other elements by the relation
\((ab)^*=b^*a^*\). Similarly the \emph{conjugation} is defined on
vectors by \(\bar{x}=-x\) and the relation
\(\overline{ab}=\bar{b}\bar{a}\). We also use the notation
\(\modulus{a}^2=a\bar{a}\) for any product \(a\) of vectors.
An important observation is that any non-zero \(x\in\Space{R}{n00}\) has a
multiplicative inverse: \(x^{-1}=\frac{\bar{x}}{\modulus{x}^2}\).
For a \(2\times 2\)-matrix \(M=
 \begin{pmatrix}
   a&b\\c&d
 \end{pmatrix}\) with Clifford entries we define,
 cf.~\cite{Cnops02a}*{(4.7)}
\begin{equation}
  \label{eq:matrix-bar-star}
  \bar{M}=
\begin{pmatrix}
  d^*&-b^*\\-c^*&a^*
\end{pmatrix}\qquad \text{ and } \qquad
M^*=\begin{pmatrix}
  \bar{d} &\bar{b}\\\bar{c}&\bar{a}
\end{pmatrix}.
\end{equation}
Then \(M\bar{M}=\delta I\) for the \emph{pseudodeterminant} \(\delta:=ad^*-bc^*\) .


Quadrics in \(\Space{R}{pq}\)---which we continue to
call cycles---can be associated to
\(2\times 2\) matrices through the FSC
construction~\citelist{\cite{FillmoreSpringer90a}
  \cite{Cnops02a}*{(4.12)} \cite{Kisil12a}*{\S~4.4}}:
\begin{equation}
  \label{eq:spheres-Rn}
  k\bar{x}x-l\bar{x}-x\bar{l}+m=0 \quad \leftrightarrow \quad
  \cycle{}{}=
  \begin{pmatrix}
    l & m\\
    k & \bar{l}
  \end{pmatrix}\notingiq
\end{equation}
where \(k, m\in\Space{R}{}\) and \(l\in\Space{R}{pq}\).  For brevity
we also encode a cycle by its coefficients \((k,l,m)\).  A
justification of~\eqref{eq:spheres-Rn} is provided by the identity:
\begin{displaymath}
  \begin{pmatrix}
    1&\bar{x}
  \end{pmatrix}
  \begin{pmatrix}
    l & m\\
    k & \bar{l}
  \end{pmatrix}
  \begin{pmatrix}
    {x}\\1
  \end{pmatrix}=
  kx\bar{x}-l\bar{x}-x\bar{l}+m,\quad \text{ since } \bar{x}=-x \text{
    for } x\in\Space{R}{pq}.
\end{displaymath}
The identification is also FLT-covariant in the sense that the
transformation~\eqref{eq:flt-defn} associated with the matrix \(M=
\begin{pmatrix}
  a&b\\c&d
\end{pmatrix}
\) sends a cycle
\(\cycle{}{}\) to the cycle \(M\cycle{}{}M^{*}\)~\cite{Cnops02a}*{(4.16)}.
We define the FLT-invariant inner product of cycles
\(\cycle{}{1}\) and \(\cycle{}{2}\) by the
identity
\begin{align}
  \label{eq:cycle-product}
  \scalar{\cycle{}{1}}{\cycle[]{}{2}}&=\Re
\tr(\cycle{}{1}\cycle{}{2})\,\notingiq
\end{align}
where \(\Re\) denotes the scalar part of a Clifford number. This
definition in term of matrices immediately implies that the inner
product is FLT-invariant. The explicit expression in terms of
 components of cycles \(\cycle{}{1}=(k_1,l_1,m_1)\) and
\(\cycle{}{2}=(k_2,l_2,m_2)\) is also useful sometimes:
\begin{align}
  \label{eq:cycle-product-expl}
  \scalar{\cycle{}{1}}{\cycle[]{}{2}}&=l_1 l_2+ \bar{l}_1 \bar{l}_2+m_1k_2+m_2k_1\,.
\end{align}
As usual, the relation \(\scalar{\cycle{}{1}}{\cycle[]{}{2}}=0\) is
called the \emph{orthogonality} of cycles \(\cycle{}{1}\) and
\(\cycle{}{2}\). In most cases it corresponds to orthogonality of
quadrics in the point space. More generally, most of
FLT-invariant relations between quadrics may be expressed in
terms FLT-invariant inner product~\eqref{eq:cycle-product}. For
the full description of methods on individual cycles, which are
implemented in the library {{\bf{}cycle}}, see the respective
documentation~\cite{Kisil05b}.

\begin{remark}
  Since cycles are elements of the projective space, the following
  \emph{normalised cycle product}:
  \begin{equation}
    \label{eq:norm-cycle-prod}
    \nscalar{C_1}{C_2}:=\frac{\scalar{C_1}{C_2}}{\sqrt{\scalar{C_1}{C_1}
        \scalar{C_2}{C_2}}}
  \end{equation}
  is more meaningful than the cycle product~\eqref{eq:cycle-product}
  itself. Note that, \(\nscalar{C_1}{C_2}\) is defined only if neither
  \(C_1\) nor \(C_2\) is a zero-radius cycle (i.e. a point). Also, the
  normalised cycle product is \(\mathrm{GL}_2(\Space{C}{})\)-invariant
  in comparison to \(\mathrm{SL}_2(\Space{C}{})\)-invariance
  of~\eqref{eq:cycle-product}.
\end{remark}

We finish this brief review of the library {{\bf{}cycle}} by pointing to
its light version written in \textsf{Asymptote}
language~\cite{Asymptote} and distributed together with the
paper~\cite{KisilReid18a}. Although the light version mostly inherited
API of the library {{\bf{}cycle}},  there are some significant limitations caused
by the absence of {\GiNaC} support:
\begin{enumerate}
\item there is no symbolic computations of any sort{\notingiqsemitocoma}
\item the light version works in two dimensions only{\notingiqsemitocoma}
\item only elliptic metrics in the point and cycle spaces are supported.
\end{enumerate}
On the other hand, being integrated with  \textsf{Asymptote} the light
version simplifies production of illustrations, which are its main target.

\subsection{Connecting quadrics and cycles}
\label{sec:conn-quadr-cycl}
The library {{\bf{}figure}} has an ability to store and resolve the system of
geometric relations between cycles. We explain below some mathematical
foundations, which greatly simplify this task.

We need a vocabulary, which translates geometric properties of
quadrics on the point space to corresponding relations in the cycle
space. The key ingredient is the cycle
product~\eqref{eq:cycle-product}--\eqref{eq:cycle-product-expl}, which
is linear in each cycles\vtick\  parameters. However, certain conditions,
e.g. tangency of cycles, involve polynomials of cycle products and
thus are non-linear.  For a successful algorithmic implementation, the
following observation is important: \emph{all non-linear conditions
  below can be linearised if the additional quadratic condition of
  normalisation type is imposed}:
\begin{equation}
  \label{eq:det-normalisation-cond}
  \scalar{\cycle{}{}}{\cycle{}{}}=\pm1.
\end{equation}
This observation in the context of the Apollonius problem was already
made in~\cite{FillmoreSpringer00a}. Conceptually the present work has
a lot in common with the above mentioned paper of Fillmore and
Springer, however a reader need to be warned that our implementation is
totally different (and, interestingly, is more closer to another
paper~\cite{FillmoreSpringer90a} of Fillmore and Springer).
\begin{remark}
  Interestingly, the method of order reduction for algebraic equations is
  conceptually similar to the method of order reduction of
  differential equations used to build a geometric dynamics of quantum
  states in~\cite{AlmalkiKisil18a}.
\end{remark}

Here is the list of relations between cycles implemented in the
current version of the library {{\bf{}figure}}.
\begin{enumerate}
\item  \label{item:quadric-flat}
  A quadric is flat (i.e. is a hyperplane), that is, its equation
  is linear. Then, either of two equivalent conditions can be used:
  \begin{enumerate}
  \item \(k\) component of the cycle vector is zero{\notingiqsemitocoma}
  \item the cycle is orthogonal
    \(\scalar{\cycle{}{1}}{\cycle[]{}{\infty}}=0\) to the ``zero-radius cycle at
    infinity'' \(\cycle[]{}{\infty}=(0,0,1)\).
  \end{enumerate}
\item \label{it:lobachevski-line}
  A quadric on the plane represents a line in Lobachevsky-type
  geometry if it is orthogonal
  \(\scalar{\cycle{}{1}}{\cycle[]{}{\Space{R}{}}}=0\)  to the real line cycle
  \(\cycle{}{\Space{R}{}}\). A similar condition is meaningful in
  higher dimensions as well.
\item \label{it:point-zero-radius}
  A quadric \(\cycle{}{}\) represents a point, that is, it has zero
  radius at given metric of the point space. Then, the determinant of
  the corresponding FSC matrix is zero or, equivalently, the cycle is
  self-orthogonal (isotropic):
  \(\scalar{\cycle{}{}}{\cycle[]{}{}}=0\). Naturally, such a cycle
  cannot be normalised to the form~\eqref{eq:det-normalisation-cond}.
\item Two quadrics are orthogonal in the point space
  \(\Space{R}{pq}\). Then, the matrices representing cycles are
  orthogonal in the sense of the inner
  product~\eqref{eq:cycle-product}.
\item Two cycles \(\cycle{}{}\) and \(\cycle[\tilde]{}{}\)
  are tangent. Then we have the following quadratic condition:
  \begin{equation}
    \label{eq:tangent-condition-defn}
    \scalar{\cycle{}{}}{\cycle[\tilde]{}{}}^2
    =  \scalar{\cycle{}{}}{\cycle{}{}}
    \scalar{\cycle[\tilde]{}{}}{\cycle[\tilde]{}{}}
    \qquad \left(\text{ or }
    \nscalar{\cycle{}{}}{\cycle[\tilde]{}{}}=\pm 1\right).
  \end{equation}
  With the assumption, that the cycle \(\cycle{}{}\) is normalised by
  the condition~\eqref{eq:det-normalisation-cond}, we may re-state
  this condition in the relation, which is linear to components of the cycle
  \(\cycle{}{}\):
  \begin{equation}
    \label{eq:tangent-condition-linear}
    \scalar{\cycle{}{}}{\cycle[\tilde]{}{}}
    = \pm \sqrt{\scalar{\cycle[\tilde]{}{}}{\cycle[\tilde]{}{}}}.
  \end{equation}
  Different signs here represent internal and outer touch.
\item Inversive distance \(\theta\) of two (non-isotropic) cycles is
  defined by the formula:
  \begin{equation}
    \label{eq:inversive-distance}
    \scalar{\cycle{}{}}{\cycle[\tilde]{}{}}
    = \theta \sqrt{
    \scalar{\cycle{}{}}{\cycle{}{}}}
  \sqrt{\scalar{\cycle[\tilde]{}{}}{\cycle[\tilde]{}{}}}
  \end{equation}
  In particular, the above discussed orthogonality corresponds to
  \(\theta=0\) and the tangency to \(\theta=\pm1\). For intersecting
  spheres \(\theta\) provides the cosine of the intersecting
  angle. For other metrics, the geometric interpretation of inversive
  distance shall be modified accordingly.

  If we are looking for a cycle \({\cycle{}{}}\) with a given
  inversive distance \(\theta\) to a given cycle
  \({\cycle[\tilde]{}{}}\), then the
  normalisation~\eqref{eq:det-normalisation-cond} again turns the
  defining relation~\eqref{eq:inversive-distance} into a linear with
  respect to parameters of the unknown cycle  \({\cycle{}{}}\).
\item A generalisation of Steiner power \(d\) of two cycles is defined
  as, cf.~\cite{FillmoreSpringer00a}*{\S~1.1}:
  \begin{equation}
    \label{eq:steiner-power}
    d=    \scalar{\cycle{}{}}{\cycle[\tilde]{}{}}
    + \sqrt{\scalar{\cycle{}{}}{\cycle{}{}}}
  \sqrt{\scalar{\cycle[\tilde]{}{}}{\cycle[\tilde]{}{}}}\notingiq
  \end{equation}
  where both cycles \(\cycle{}{}\) and \(\cycle[\tilde]{}{}\) are
  \(k\)-normalised, that is the coefficient in front the quadratic
  term in~\eqref{eq:spheres-Rn} is \(1\). Geometrically, the
  generalised Steiner power for spheres provides the square of
  tangential distance. However, this relation is again non-linear for
  the cycle \(\cycle{}{}\).

  If we replace \(\cycle{}{}\) by the cycle
  \(\cycle{}{1}=\frac{1}{\sqrt{\scalar{\cycle{}{}}{\cycle{}{}}}}\cycle{}{}\)
  satisfying~\eqref{eq:det-normalisation-cond}, the
  identity~\eqref{eq:steiner-power} becomes:
  \begin{equation}
    \label{eq:steiner-power-linear}
    d\cdot k=    \scalar{\cycle{}{1}}{\cycle[\tilde]{}{}}
    +   \sqrt{\scalar{\cycle[\tilde]{}{}}{\cycle[\tilde]{}{}}}\notingiq
  \end{equation}
  where \(k=\frac{1}{\sqrt{\scalar{\cycle{}{}}{\cycle{}{}}}}\) is the
  coefficient in front of the quadratic term of \(\cycle{}{1}\). The
  last identity is linear in terms of the coefficients of
  \(\cycle{}{1}\).
\end{enumerate}
Summing up: if an unknown cycle is connected to already given cycles
by any combination of the above relations, then all conditions can be
expressed as \emph{a system of linear equations for coefficients of the
unknown cycle and at most one quadratic
equation~\eqref{eq:det-normalisation-cond}}.

\section[Example: Poincare Extension of Moebius Transformations]{Example: Poincar\'e Extension of M\"obius Transformations}
\label{sec:poinc-extens-mobi}

It is known, that M\"obius transformations on \(\Space{R}{n}\) can be
expanded to the ``upper'' half-space in \(\Space{R}{n+1}\) using the
Poincar\'e extension~\citelist{\cite{Beardon95}*{\S~3.3}
  \cite{JParker07a}*{\S~5.2}}.  An explicit formula is usually
presented without a discussion of its origin. In particular, one may
get an impression that the solution is
unique. Following~\cite{Kisil15a} we consider various aspects of such
extension and describe different possible realisations. Our
consideration is restricted to the case of extension from the real
line to the upper half-plane. However, we made an effort to present it
in a way, which allows numerous further generalisations. 
In section~\ref{sec:figures-as-families} we a partial realisation of
Poincar\'e extension will be formulated through ensembles of cycles.

\subsection{Geometric construction}
\label{sec:geom-constr}

We start from the geometric procedure in the standard
situation. Recall, the
group \(\SL\) consists of  real \(2\times 2\) matrices with the unit determinant.
\(\SL\) acts on the real line by linear-fractional maps:
\begin{equation}
  \label{eq:moebius-map-defn}
  \begin{pmatrix}
    a&b\\c&d
  \end{pmatrix}: x \mapsto \frac{ax+b}{cx+d},
  \quad \text{ where } 
  x\in\Space{R}{}\text{ and } 
    \begin{pmatrix}
    a&b\\c&d
  \end{pmatrix}\in\SL.
\end{equation}
A pair of (possibly equal) real numbers \(x\)
and \(y\)
uniquely determines a semicircle \(\cycle{}{xy}\)
in the upper half-plane with the diameter \([x,y]\).
For a linear-fractional transformation
\(M\)~\eqref{eq:moebius-map-defn},
the images \(M(x)\)
and \(M(y)\)
define the semicircle with the diameter \([M(x),M(y)]\),
thus, we can define the action of \(M\)
on semicircles: \(M(\cycle{}{xy})=\cycle{}{M(x)M(y)}\).
Geometrically, the Poincar\'e extension is based on the following
lemma, see Fig.~\ref{fig:poincare-extension}(a) and more general
Lem.~\ref{le:poincare-geom-general} below:
\begin{lemma}
  \label{le:pencil-common-point}
  If a pencil of semicircles in the upper half-plane has a common
  point, then the images of these semicircles under a
  transformation~\eqref{eq:moebius-map-defn} have a common point as
  well.
\end{lemma}
Elementary geometry of right triangles tells that a pair of
intersecting intervals \([x,y]\), \([x',y']\), where \(x<x'<y<y'\),
defines the point
\begin{equation}
  \label{eq:poincare-point-ell}
  \left( \frac{xy-x'y'}{x+y-x'-y'}\,,\,
    \frac{\sqrt{(x-y')(x-x')(x'-y)(y-y')}}{x+y-x'-y'}\right) \in \Space[+]{R}{2}.
\end{equation}
An alternative demonstration uses three observations:
\begin{enumerate}
\item the scaling \(x\mapsto ax\),
  \(a>0\)
  on the real line produces the scaling \((u,v) \mapsto(au,av)\)
  on pairs~\eqref{eq:poincare-point-ell}{\notingiqsemitocoma}
\item the horizontal shift \(x\mapsto x+b\) on the real line produces
  the horizontal shift \((u,v) \mapsto (u+b,v)\) on
  pairs~\eqref{eq:poincare-point-ell}{\notingiqsemitocoma} 
\item for the special case \(y=-x^{-1}\) and \(y'=-{x'}^{-1}\) 
  the pair~\eqref{eq:poincare-point-ell} is \((0,1)\). 
\end{enumerate}
Finally, expression~\eqref{eq:poincare-point-ell}, as well
as~\eqref{eq:poincare-point-hyp}--\eqref{eq:poincare-point-par} below,
can be calculated by the specialised CAS for M\"obius invariant
geometries~\cites{Kisil05b,Kisil14b}.

\begin{figure}[htbp]
  \centering
  \makebox[0pt][l]{(a)}\,\includegraphics[scale=.85]{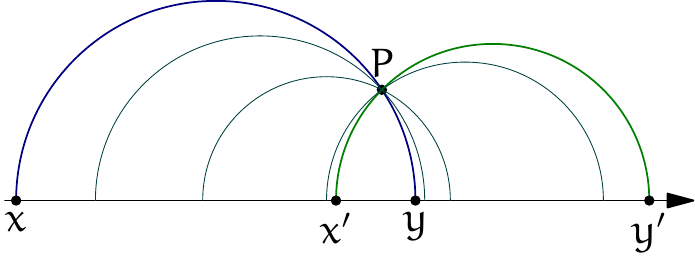}\hfill
  \makebox[0pt][l]{(b)}\,\includegraphics[scale=.85]{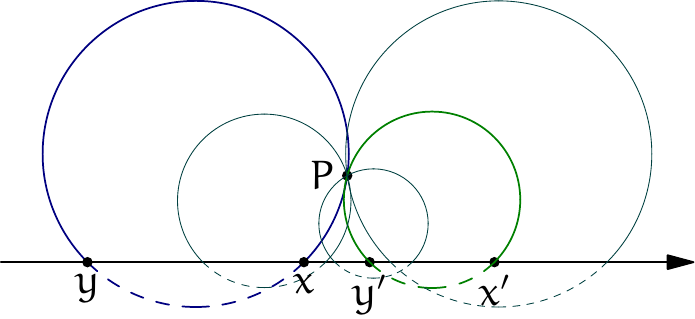}\\[2em]
  \makebox[0pt][l]{(c)}\,\includegraphics[scale=.85]{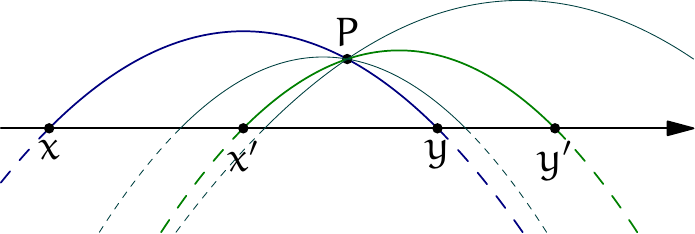}\hfill
  \makebox[0pt][l]{(d)}\,\includegraphics[scale=.85]{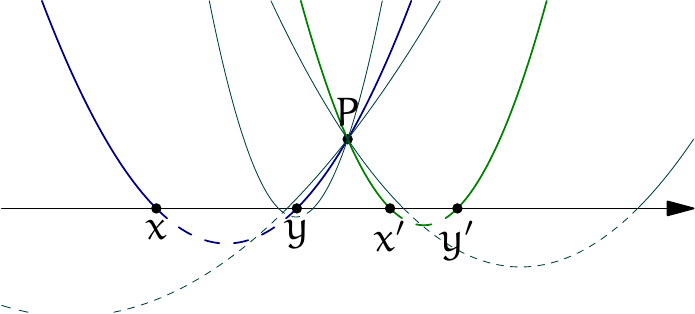}\\[2em]
  \makebox[0pt][l]{(e)}\,\includegraphics[scale=.85]{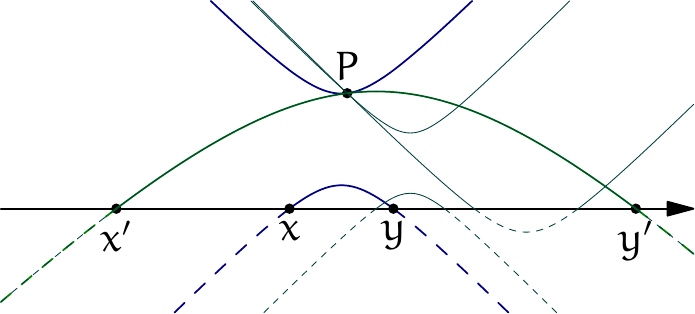}\hfill
  \makebox[0pt][l]{(f)}\,\includegraphics[scale=.85]{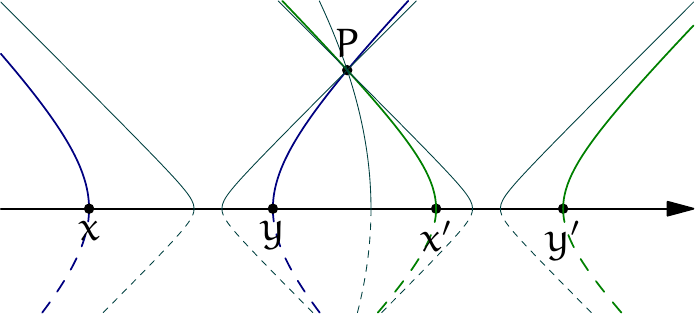}
  \caption[Poincar\'e extensions]{Poincar\'e extensions: first column
    presents points defined by the intersecting intervals \([x,y]\)
    and \([x',y']\), the second column---by disjoint intervals. Each
    row uses the same type of conic sections---circles, parabolas and
    hyperbolas respectively. Pictures are produced by the software~\cite{Kisil14b}.}
  \label{fig:poincare-extension}
\end{figure}

This standard approach can be widened as follows. The above semicircle
can be equivalently described through the unique circle passing \(x\)
and \(y\)
and orthogonal to the real axis. Similarly, an interval \([x,y]\)
uniquely defines a right-angle hyperbola in \(\Space{R}{2}\)
orthogonal to the real line and passing (actually, having her vertices
at) \((x,0)\)
and \((y,0)\).
An intersection with the second such hyperbola having vertices \((x',0)\)
and \((y',0)\)
defines a point with coordinates, see
Fig.~\ref{fig:poincare-extension}(f):
\begin{equation}
  \label{eq:poincare-point-hyp}
  \left( \frac{xy-x'y'}{x+y-x'-y'}\,,\,
    \frac{\sqrt{(x-y')(x-x')(x'-y)(y'-y)}}{x+y-x'-y'}\right), 
\end{equation}
where \(x<y<x'<y'\).  Note, the opposite sign of the product under the
square roots in~\eqref{eq:poincare-point-ell}
and~\eqref{eq:poincare-point-hyp}.

If we wish to consider the third type of conic
sections---parabolas---we cannot use the unaltered procedure: there is
no a non-degenerate parabola orthogonal to the real line and
intersecting the real line in two points. We may recall, that a circle
(or hyperbola) is orthogonal to the real line if its centre belongs to
the real line. Analogously, a parabola is \emph{focally orthogonal} (see
\cite{Kisil12a}*{\S~6.6} for a general consideration) to the real line
if its focus belongs to the real line. Then, an
interval \([x,y]\) uniquely defines a downward-opened parabola with
the real roots \(x\) and \(y\) and focally orthogonal to the real line.  Two
such parabolas defined by intervals \([x,y]\) and \([x',y']\) have a
common point, see Fig.~\ref{fig:poincare-extension}(c):
\begin{equation}
  \label{eq:poincare-point-par}
  \left(
    \frac{xy'-yx'+D}{x-y-x'+y'}\,,\,
     \frac{(x'-x) (y'-y) (y-x+y'-x')+(x+y-x'-y') D}{(x-y-x'+y')^{2}}
\right)\notingiq
\end{equation}
where \(D=\pm\sqrt{(x-x')(y-y')(y-x)(y'-x')}\). 
For  pencils of such hyperbolas and parabolas respective variants of
Lem.~\ref{le:pencil-common-point} hold.

Focally orthogonal parabolas make the angle \(45^\circ\) with the real
line. This suggests to replace orthogonal circles and hyperbolas by
conic sections with a fixed angle to the real line, see
Fig.~\ref{fig:poincare-extension}(b)--(e). Of course, to be consistent
this procedure requires a suitable modification of
Lem.~\ref{le:pencil-common-point}, we will obtain it as a byproduct
of our study, see Lem.~\ref{le:poincare-geom-general}. However, the
respective alterations of the above formulae
\eqref{eq:poincare-point-ell}--\eqref{eq:poincare-point-par} become
more complicated in the general case.


The considered geometric construction is elementary and visually
appealing. Now we turn to respective algebraic consideration. 

\subsection[Moebius transformations and Cycles]{M\"obius transformations and Cycles}
\label{sec:cycles}
The group \(\SL\) acts on \(\Space{R}{2}\) by matrix multiplication on
column vectors:
\begin{equation}
  \label{eq:left-matrix-mult}
  \oper{L}_g:\ 
  \begin{pmatrix}
    x_1\\x_2
  \end{pmatrix}
  \mapsto
  \begin{pmatrix}
    ax_1+bx_2\\cx_1+dx_2
  \end{pmatrix}=
  \begin{pmatrix}
    a&b\\c&d
  \end{pmatrix}  \begin{pmatrix}
    x_1\\x_2
  \end{pmatrix}
  \,, \quad \text{ where } g=\begin{pmatrix}
    a&b\\c&d
  \end{pmatrix}\in\SL.
\end{equation}
A linear action respects the equivalence relation \(  \begin{pmatrix}
  x_1\\x_2
\end{pmatrix}
\sim   \begin{pmatrix}
  \lambda x_1\\ \lambda x_2
\end{pmatrix}\), \(\lambda\neq 0\) on \(\Space{R}{2}\). The
collection of all cosets for non-zero vectors in \(\Space{R}{2}\) is the \emph{projective
  line} \(P\Space{R}{1}\). Explicitly, a
non-zero vector \(\begin{pmatrix} x_1\\x_2
\end{pmatrix}\in\Space{R}{2}\) corresponds to the point with
homogeneous coordinates
\([x_1:x_2]\in P\Space{R}{1}\). If \(x_2\neq 0\) then this point is
represented by \([\frac{x_1}{x_2}:1]\) as well. 
The embedding \(\Space{R}{} \rightarrow P\Space{R}{1}\) defined by  \(x \mapsto
[x:1]\), 
\(x\in\Space{R}{}\) covers the all but one
of points in \(P\Space{R}{1}\). The exceptional point \([1:0]\) is
naturally identified with the infinity.

The linear action~\eqref{eq:left-matrix-mult} induces a morphism of
the projective line \(P\Space{R}{1}\),
which is called a M\"obius transformation. Considered on the real line
within \(P\Space{R}{1}\),
M\"obius transformations takes fraction-linear form:
\begin{displaymath}
  g:\ 
  [x:1]
  \mapsto
  \left[\frac{ax+b}{cx+d}:1\right]\,, \quad \text{ where } g=\begin{pmatrix}
    a&b\\c&d
  \end{pmatrix}\in\SL \text{ and } cx+d\neq 0.
\end{displaymath}
This \(\SL\)-action on \(P\Space{R}{1}\) is denoted as \(g: x\mapsto
g\cdot x\).  We note that the correspondence of
column vectors and row vectors \( 
\map{i}: \begin{pmatrix}
  x_1\\x_2
\end{pmatrix} \mapsto (x_2, -x_1)\) intertwines
the left multiplication \(\oper{L}_g\)~\eqref{eq:left-matrix-mult} and
the right multiplication \(\oper{R}_{g^{-1}}\) by the inverse matrix:
\begin{equation}
  \label{eq:right-matrix-mult}
  \oper{R}_{g^{-1}}:\ 
    (x_2,-x_1)
  \mapsto
    (cx_1+dx_2,\, -ax_1-bx_2)
    =
    (x_2,-x_1)
  \begin{pmatrix}
    d&-b\\-c&a
  \end{pmatrix}\,.
\end{equation}
We extended the map \(\map{i}\) to \(2\times 2\)-matrices by the rule:
\begin{equation}
  \label{eq:t-map-matrix}
  \map{i}:\
  \begin{pmatrix}
    x_1&y_1\\
    x_2&y_2
  \end{pmatrix}
  \mapsto
  \begin{pmatrix}
    y_2&-y_1\\
    x_2&-x_1
  \end{pmatrix}\,.
\end{equation}
Two columns \(\begin{pmatrix} x\\1
\end{pmatrix}\) and \(\begin{pmatrix} y\\1
\end{pmatrix}\) form the \(2\times 2\) matrix \(M_{xy}=\begin{pmatrix} x&y\\1&1
\end{pmatrix}\). For geometrical reasons appearing in
Cor.~\ref{co:null-set-quadratic}, we call a \emph{cycle} the
\(2\times 2\)-matrix \(\cycle{}{xy}\) defined by
\begin{equation}
  \label{eq:matrix-for-cycle-defn}
\cycle{}{xy}=\frac{1}{2} M_{xy}\cdot\map{i}(M_{xy})=\frac{1}{2} M_{yx}\cdot\map{i}(M_{yx})=
  \begin{pmatrix}
    \frac{x+y}{2}  & -xy\\
    1 & -\frac{x+y}{2}
  \end{pmatrix}\,.
\end{equation}
We note that \(\det \cycle{}{xy}=-(x-y)^2/4\), thus \(\det \cycle{}{xy}=0\) if and
only if \(x=y\). Also, we can consider the M\"obius transformation
produced by the \(2\times 2\)-matrix \(\cycle{}{xy}\) and calculate:
\begin{equation}
  \label{eq:eigenvectors-c_xy}
\cycle{}{xy}
  \begin{pmatrix}
  x\\1  
  \end{pmatrix}=
\lambda   \begin{pmatrix}
  x\\1  
  \end{pmatrix}
\quad \text{and } \quad
\cycle{}{xy}
  \begin{pmatrix}
  y\\1  
  \end{pmatrix}=
-\lambda   \begin{pmatrix}
  y\\1  
  \end{pmatrix}
\quad\text{where } \lambda=\frac{x-y}{2}.
\end{equation}
Thus, points \([x:1]\), \([y:1]\in P\Space{R}{1}\) are fixed by
\(\cycle{}{xy}\). Also,  \(\cycle{}{xy}\) swaps the interval \([x,y]\) and its
complement. 

Due to their structure, matrices \(\cycle{}{xy}\) can be parametrised by
points of \(\Space{R}{3}\).  Furthermore, the map from
\(\Space{R}{2}\rightarrow \Space{R}{3}\) given by \((x,y)\mapsto
\cycle{}{xy}\) naturally induces the projective map \((P\Space{R}{1})^2
\rightarrow P\Space{R}{2}\) due to the identity:
\begin{displaymath}
\frac{1}{2}  \begin{pmatrix}
    \lambda x & \mu y\\
    \lambda & \mu
  \end{pmatrix}
  \begin{pmatrix}
    \mu & -\mu y\\
    \lambda & - \lambda x
  \end{pmatrix}=\lambda \mu    \begin{pmatrix}
    \frac{x+y}{2} & -xy\\
    1 & -\frac{x+y}{2}
  \end{pmatrix}= \lambda \mu \cycle{}{xy}\,.
\end{displaymath}
Conversely, a zero-trace matrix \(
\begin{pmatrix}
  a&b\\c&-a
\end{pmatrix}
\) with a non-positive determinant is projectively equivalent to
a product \(\cycle{}{xy}\)~\eqref{eq:matrix-for-cycle-defn} with
\(x,y=\frac{a\pm\sqrt{a^2+bc}}{c}\). In particular, we can embed a point
 \([x:1]\in P\Space{R}{1}\) to \(2\times 2\)-matrix \(\cycle{}{xx}\) with
 zero determinant. 

The combination
of~\eqref{eq:left-matrix-mult}--\eqref{eq:matrix-for-cycle-defn} implies
that the correspondence \((x,y)\mapsto \cycle{}{xy}\) is \(\SL\)-covariant in the
following sense:
\begin{equation}
  \label{eq:covariance-Cxy}
  g \cycle{}{xy}  g^{-1}=\cycle{}{x'y'}\,,\quad \text{ where }
  x'=g\cdot x \text{ and } y'=g\cdot y.
\end{equation}

To achieve a geometric interpretation of all matrices, we consider the
bilinear form \(Q: \Space{R}{2}\times \Space{R}{2}\rightarrow
\Space{R}{}\) generated by a \(2\times 2\)-matrix \(\begin{pmatrix}
    a&b\\
    c&d
  \end{pmatrix}\):
\begin{equation}
  \label{eq:bilinear-form-defn}
  Q(x, y)=
  \begin{pmatrix}
    x_1 & x_2
  \end{pmatrix}
  \begin{pmatrix}
    a&b\\
    c&d
  \end{pmatrix}
  \begin{pmatrix}
    y_1 \\ y_2
  \end{pmatrix}\, , \quad \text{ where } x=(x_1,x_2),\ y=(y_1,y_2). 
\end{equation}
Due to linearity of \(Q\), the null set 
\begin{equation}
  \label{eq:bilin-form-null-set}
  \{(x,y)\in \Space{R}{2}\times \Space{R}{2} \such Q(x,y)=0\}
\end{equation}
factors to a subset of \(P\Space{R}{1}\times
P\Space{R}{1}\). Furthermore, for the matrices
\(\cycle{}{xy}\)~\eqref{eq:matrix-for-cycle-defn}, a direct calculation
shows that:
\begin{lemma}
  \label{lem:biliniar-inner}
  The following identity holds:
  \begin{equation}
    \label{eq:cycle-inner-product}
    \cycle{}{xy}(\map{i}(x'),y')=\tr (\cycle{}{xy} \cycle{}{x'y'})
      =\textstyle  \frac{1}{2}(x+y)(x'+ y')-(x y+x' y')\,.
  \end{equation}
  In particular, the above expression is a symmetric function of the pairs
  \((x,y)\) and \((x',y')\).
\end{lemma}
The map \(\map{i}\) appearance in~\eqref{eq:cycle-inner-product} is
justified once more by the following result.
\begin{corollary}
  \label{co:null-set-quadratic}
  The null set of the quadratic form
  \(\cycle{}{xy}(x')=\cycle{}{xy}(\map{i}(x'),x')\) consists of two points \(x\)
  and \(y\).
\end{corollary}
Alternatively, the identities \(\cycle{}{xy}(x)=\cycle{}{xy}(y)=0\) follows
from~\eqref{eq:eigenvectors-c_xy} and the fact that \(\map{i}(z)\) is
orthogonal to \(z\) for all \(z\in\Space{R}{2}\).  Also, we note that:
\begin{displaymath}
  \map{i}\left(\!\begin{pmatrix}
    x_1 \\ x_2
  \end{pmatrix}\!\right)
  \begin{pmatrix}
    a&b\\
    c&d
  \end{pmatrix}
  \begin{pmatrix}
    y_1 \\ y_2
  \end{pmatrix}=
  \begin{pmatrix}
    x_1 & x_2
  \end{pmatrix}
  \begin{pmatrix}
    -c&-d\\
    a&b
  \end{pmatrix}
  \begin{pmatrix}
    y_1 \\ y_2
  \end{pmatrix}.
\end{displaymath}

Motivated by Lem.~\ref{lem:biliniar-inner}, we call
\(\scalar{\cycle{}{xy}}{\cycle{}{x'y'}}:=-\tr (\cycle{}{xy} \cycle{}{x'y'})\) the
\emph{pairing} of two cycles.  It shall be noted that
the pairing is \emph{not} positively defined, this follows from the
explicit expression~\eqref{eq:cycle-inner-product}. The sign is chosen
in such way, that
\begin{displaymath}
\scalar{\cycle{}{xy}}{\cycle{}{xy}}=-2\det(\cycle{}{xy})=\textstyle \frac{1}{2}(x-y)^2\geq 0.  
\end{displaymath}
Also, an immediate consequence of Lem.~\ref{lem:biliniar-inner} or
identity~\eqref{eq:bilinear-form-defn} is
\begin{corollary}
  The  pairing of cycles is invariant under the 
  action~\eqref{eq:covariance-Cxy} of \(\SL\):
  \begin{displaymath}
    \scalar{g\cdot \cycle{}{xy} \cdot g^{-1} }{g\cdot \cycle{}{x'y'} \cdot g^{-1} } =\scalar{\cycle{}{xy}}{\cycle{}{x'y'}}.
  \end{displaymath}
\end{corollary}
From~\eqref{eq:cycle-inner-product}, the null
set~\eqref{eq:bilin-form-null-set} of the form \(Q=\cycle{}{xy}\) can be
associated to the family of cycles \(\{\cycle{}{x'y'} \such
\scalar{\cycle{}{xy}}{\cycle{}{x'y'}}=0, (x',y')\in\Space{R}{2}\times
\Space{R}{2}\}\) which we will call \emph{orthogonal} to \(\cycle{}{xy}\).

\subsection{Extending cycles}
\label{sec:extending-cycles}

Since bilinear forms with matrices \(\cycle{}{xy}\) have numerous geometric
connections with \(P\Space{R}{1}\), we are looking for a similar
interpretation of the generic matrices. The previous discussion
identified the key ingredient of the recipe: \(\SL\)-invariant
pairing~\eqref{eq:cycle-inner-product} of two forms. Keeping in mind
the structure of \(\cycle{}{xy}\), we will parameterise\footnote{Further
  justification of this parametrisation will be obtained from the
  equation of a quadratic curve~\eqref{eq:quadratics-equation}.} a
generic \(2\times 2\) matrix as \(\begin{pmatrix} l+n&-m\\k&-l+n
\end{pmatrix}\) and consider the corresponding four dimensional vector
\((n,l,k,m)\). Then, the similarity with \(\begin{pmatrix}
    a&b\\c&d
  \end{pmatrix}
  \in\SL\):
\begin{equation}
  \label{eq:sl2-similarity-2x2-matrices}
  \begin{pmatrix}
    {l}'+{n}'&-{m}'\\{k}'&-{l}'+{n}'
  \end{pmatrix}
  =
  \begin{pmatrix}
    a&b\\c&d
  \end{pmatrix}
  \begin{pmatrix}
    l+n&-m\\k&-l+n
  \end{pmatrix}
  \begin{pmatrix}
    a&b\\c&d
  \end{pmatrix}^{-1}\notingiq
\end{equation}
corresponds to the linear transformation of \(\Space{R}{4}\),
cf.~\cite{Kisil12a}*{Ex.~4.15}: 
\begin{equation}
  \label{eq:SL2-act-cycle-linear}
  \begin{pmatrix}
    {n}'\\{l}'\\{k}'\\{m}'
  \end{pmatrix}=
  \begin{pmatrix}
    1&0&0&0\\
    0&c b+a d& b d&c a\\
    0&2 c d&d^2&c^2\\
    0&2 a b&b^2&a^2
  \end{pmatrix}
  \begin{pmatrix}
    n\\l\\k\\m
  \end{pmatrix}.
\end{equation}
In particular, this action commutes with the scaling of the first component:
\begin{equation}
  \label{eq:scaling-n}
  \lambda: (n,l,k,m) \mapsto  (\lambda n,l,k,m).
\end{equation}
This expression is helpful in proving the following statement.
\begin{lemma}
  \label{le:sl2-invariant-pairing}
  Any \(\SL\)-invariant (in the sense of the
  action~\eqref{eq:SL2-act-cycle-linear}) pairing in \(\Space{R}{4}\) is
  isomorphic to
  \begin{displaymath}
  2{\bsigma} n n'-2{l}l'+{k}m'+{m}k'
  =
  \begin{pmatrix}
    {n}'&{l}'&{k}'&{m}'
  \end{pmatrix}
  \begin{pmatrix}
    2{\bsigma}&0&0&0\\0&-2&0&0\\0&0&0&1\\0&0&1&0
  \end{pmatrix}
  \begin{pmatrix}
    n\\l\\k\\m
  \end{pmatrix}\notingiq
  \end{displaymath}
  where \(\bsigma =-1\), \(0\) or \(1\) and 
  \((n,l,k,m)\),  \((n',l',k',m')\in\Space{R}{4}\) .
\end{lemma}
\begin{proof}
  Let \(T\) be \(4\times 4\) a matrix
  from~\eqref{eq:SL2-act-cycle-linear}, if a \(\SL\)-invariant pairing
  is defined by a \(4\times 4\) matrix \(J=(j_{fg})\), then
  \(T'JT=J\), where \(T'\) is transpose of \(T\). The equivalent
  identity \(T'J=JT^{-1}\) produces a system of homogeneous linear
  equations which has the generic solution:
  \begin{align*}
    j_{12}&=
    j_{13}=
    j_{14}=
    j_{21}=
    j_{31}=
    j_{41}=0,\\
    j_{22}&= \frac{ (d-a) j_{42}-2 b j_{44}}{c}-2  j_{43},&
    j_{23}&=- \frac{b}{c} j_{42},&
    j_{24}&=-\frac{ c j_{42}+2 (a-d) j_{44}}{c},\\
    j_{34}&=\frac{c (a - d) j_{42}+(a -d)^2j_{44}}{ c^2}+ j_{43},&
    j_{33}&=  \frac{b^2}{c^2} j_{44} ,&
    j_{32}&= \frac{b}{c}\frac{c  j_{42}+2  (a-d) j_{44}}{c}\notingiq
  \end{align*}
  with four free variables \(j_{11}\), \(j_{42}\), \(j_{43}\) and
  \(j_{44}\). Since a solution shall not depend on \(a\), \(b\),
  \(c\), \(d\), we have to put \(j_{42}=j_{44}=0\). Then by the
  homogeneity of the identity \(T'J=JT^{-1}\), we can scale \(j_{43}\)
  to \(1\). Thereafter,  an independent (sign-preserving)
  scaling~\eqref{eq:scaling-n} of \(n\) leaves only three
  non-isomorphic values \(-1\), \(0\), \(1\) of \(j_{11}\).
\end{proof}
The appearance of the three essential different cases \(\bsigma =-1\),
\(0\) or \(1\) in Lem.~\ref{le:sl2-invariant-pairing} is a
manifestation of the common division of mathematical objects into
elliptic, parabolic and hyperbolic cases \citelist{\cite{Kisil06a}
  \cite{Kisil12a}*{Ch.~1}}.  Thus, we will use letters ``e'', ``p'',
``h'' to encode the corresponding three values of \(\bsigma\). 

Now we may describe all \(\SL\)-invariant pairings of bilinear forms.
\begin{corollary}
  Any \(\SL\)-invariant (in the sense of the
  similarity~\eqref{eq:sl2-similarity-2x2-matrices}) pairing between
  two bilinear forms \(Q=  \begin{pmatrix}
    l+n&-m\\k&-l+n
  \end{pmatrix}  \) and \(Q'=  \begin{pmatrix}
    l'+n'&-m'\\k'&-l'+n'
  \end{pmatrix}  \) is isomorphic to:
  \begin{align}
    \label{eq:cycle-product-expl1}
    \scalar[\tau]{Q}{Q'}&=-\tr ({Q}_\tau Q')\\
    & =
    2\tau n'n-2l'l+k'm+m'k, \quad\text{ where }
    {Q}_\tau=  \begin{pmatrix}
      {l}-\tau{n}&-{m}\\{k}&-{l}-\tau{n}, \nonumber 
    \end{pmatrix}
  \end{align}
  and \(\tau=-1\), \(0\) or \(1\).
\end{corollary}
Note that we can explicitly write \(Q_\tau\) for \(Q=
\begin{pmatrix}
  a&b\\c&d
\end{pmatrix}\) as follows:
\begin{displaymath}
  Q_e=\begin{pmatrix}
    a&b\\c&d
  \end{pmatrix},\quad
  Q_p=\begin{pmatrix}
    \frac{1}{2}(a-d)&b\\c&-\frac{1}{2}(a-d)
  \end{pmatrix},\quad
  Q_h=\begin{pmatrix}
    -d&b\\c&-a
  \end{pmatrix}.
\end{displaymath}
In particular, \(Q_h=-Q^{-1}\) and
\(Q_p=\frac{1}{2}(Q_e+Q_h)\). Furthermore, \(Q_p\) has the same
structure as \(\cycle{}{xy}\).  Now, we are ready to extend the
projective line \(P\Space{R}{1}\) to two dimensions using the analogy
with properties of cycles \(\cycle{}{xy}\).
\begin{definition}
  \label{de:product-cycles}
  \begin{enumerate}
  \item Two bilinear forms \(Q\) and \(Q'\) are
    \(\tau\)-\emph{orthogonal} if \(\scalar[\tau]{Q}{Q'}=0\). 
  \item \label{it:isotropic}
    A form is \(\tau\)-\emph{isotropic} if it is
    \(\tau\)-orthogonal to itself.
  \end{enumerate}
\end{definition}
If a form \( Q= \begin{pmatrix}
    l+n&-m\\k&-l+n
  \end{pmatrix}
\) has \(k\neq 0\) then we can scale it  to obtain \(k=1\), this form
of \(Q\)
is called \emph{normalised}. A normalised \(\tau\)-isotropic form is
completely determined by its diagonal values: \(\begin{pmatrix}
  u+v&-u^2+\tau v^2\\1&-u+v
\end{pmatrix}\). Thus, the set of such forms is in a one-to-one
correspondence with points of \(\Space{R}{2}\). Finally, a form \(
Q= \begin{pmatrix} l+n&-m\\k&-l+n
  \end{pmatrix}\) is e-orthogonal to the \(\tau\)-isotropic form \(\begin{pmatrix}
  u+v&-u^2+\tau v^2\\1&-u+v
\end{pmatrix}\) if:
\begin{equation}
  \label{eq:quadratics-equation}
  k(u^2-\tau v^2)-2lu-2nv+m=0\notingiq
\end{equation}
that is the point \((u,v)\in\Space{R}{2}\) belongs to the quadratic curve
with coefficients \((k,l,n,m)\).

\subsection{Homogeneous spaces of cycles}
\label{sec:homog-spac-cycl}
Obviously, the group \(\SL\) acts on \(P\Space{R}{1}\) transitively,
in fact it is even \(3\)-transitive in the following sense. We say
that a triple \(\{x_1,x_2,x_3\}\subset P\Space{R}{1}\) of distinct
points is \emph{positively oriented} if
\begin{displaymath}
  \text{ either }\quad x_1<x_2<x_3, \quad\text{ or }\quad  x_3<x_1< x_2\notingiq
\end{displaymath}
where we agree that the ideal point \(\infty\in P\Space{R}{1}\) is
greater than  any \(x\in\Space{R}{}\). Equivalently, 
a triple \(\{x_1,x_2,x_3\}\) of reals is positively oriented if:
\begin{displaymath}
  (x_1-x_2)(x_2-x_3)(x_3-x_1)>0.
\end{displaymath}
Also, a triple of distinct points, which is not positively oriented, is
\emph{negatively oriented}. 
A simple calculation based on the resolvent-type identity:
\begin{displaymath}
  \frac{ax+b}{cx+d}-\frac{ay+b}{cy+d}=\frac{(x-y)(ad-bc)}{(cx+b)(cy+d)}
\end{displaymath}
shows that both the positive and negative orientations of triples are
\(\SL\)-invariant. On the other hand, the reflection \(x\mapsto -x\)
swaps orientations of triples. Note, that the reflection is a M\"obius
transformation associated to the cycle
\begin{equation}
  \label{eq:matrix-J-defn}
 \cycle{}{0\infty}=\begin{pmatrix}
  1&0\\0&-1
\end{pmatrix}, \quad \text{ with } \det \cycle{}{0\infty}=-1.
\end{equation}

A significant amount of information about M\"obius transformations
follows from the fact, that any continuous one-parametric subgroup of
\(\SL\) is conjugated to one of the three following
subgroups\footnote{A reader may know that \(A\), \(N\) and \(K\) are
  factors in the Iwasawa decomposition \(\SL=ANK\)
  (cf. Cor.~\ref{co:Iwasawa-decomp}), however this important result
  does not play any r\^ole in our consideration.}:
\begin{equation}
  \label{eq:ank-subgroup}
  A=\left\{\begin{pmatrix} e^{-t} & 0\\0&e^{t}
    \end{pmatrix}\right\},\quad
  N=\left\{\begin{pmatrix} 1&t \\0&1
    \end{pmatrix}\right\},\quad
  K=\left\{ \begin{pmatrix}
      \cos t &  -\sin t\\
      \sin t & \cos t
    \end{pmatrix}\right\}\notingiq
\end{equation}
where \(t\in\Space{R}{}\). Also, it is useful to introduce
subgroups \(\Aprime\) and \(N'\) conjugated to \(A\) and \(N\)
respectively:
\begin{equation}
  \label{eq:a1n1-subgroup}
  \Aprime=\left\{\begin{pmatrix} \cosh t & \sinh t\\ \sinh t & \cosh t 
    \end{pmatrix} \such t\in \Space{R}{}\right\},\qquad
  N'=\left\{\begin{pmatrix} 1&0 \\t&1
    \end{pmatrix} \such t\in \Space{R}{}\right\}.
\end{equation}
Thereafter, all three one-parameter subgroups \(\Aprime\), \(N'\)
and \(K\) consist of all matrices with the universal structure
\begin{equation}
  \label{eq:a1n1k-universal}
\begin{pmatrix}
  a&\tau b\\b &a
\end{pmatrix}
\quad  \text{where }\tau =1,\, 0,\, -1 \text{ for }\Aprime,\, N'\,
\text{ and }
K \text{ respectively}.
\end{equation}
We use the notation \(H_\tau\) for these subgroups. Again, any
continuous one-di\-men\-sio\-nal
subgroup of \(\SL\) is conjugated to
\(H_\tau\) for an appropriate \(\tau\).

We note, that matrices from \(A\), \(N\) and \(K\) with \(t\neq 0\)
have two, one and none different real eigenvalues
respectively. Eigenvectors in \(\Space{R}{2}\) correspond to fixed
points of M\"obius transformations on \(P\Space{R}{1}\). Clearly, the
number of eigenvectors (and thus fixed points) is limited by the
dimensionality of the space, that is two. For this reason, if \(g_1\)
and \(g_2\) take equal values on three different points of
\(P\Space{R}{1}\), then \(g_1=g_2\).

Also, eigenvectors provide an effective classification tool:
\(g\in\SL\) belongs to a one-dimensional continuous subgroup conjugated
to \(A\), \(N\) or \(K\) if and only if the characteristic polynomial
\(\det(g-\lambda I)\) has two, one and none different real root(s)
respectively. We will illustrate an application of fixed points
techniques through the following well-known result, which will be used
later.
\begin{lemma}
  \label{le:three-transitive}
  Let \(\{x_1,x_2,x_3\}\) and \(\{y_1,y_2,y_3\}\) be positively
  oriented triples of points in \(\dSpace{R}{}\). Then, there is a
  unique (computable!) M\"obius map \(\phi\in\SL\) with \(\phi(x_j)=y_j\) for
  \(j=1\), \(2\), \(3\).
\end{lemma}
\begin{proof}
  Often, the statement is quickly demonstrated through an explicit
  expression for \(\phi\), cf.~\cite{Beardon05a}*{Thm.~13.2.1}. We
  will use properties of the subgroups \(A\), \(N\) and \(K\) to
  describe an algorithm to find such a map.  First, we notice that it
  is sufficient to show the Lemma for the particular case \(y_1=0\),
  \(y_2=1\), \(y_3=\infty\). The general case can be obtained from
  composition of two such maps. Another useful observation is that the
  fixed point for \(N\), that is \(\infty\), is also a fixed point of
  \(A\).

  Now, we will use subgroups \(K\), \(N\) and \(A\) in order of
  increasing number of their fixed points. First, for any \(x_3\) the
  matrix \(g'=\begin{pmatrix}
    \cos t &  \sin t\\
    -\sin t & \cos t
  \end{pmatrix}\in K\)  such that \(\cot t=x_3\)  maps \(x_3\) to
  \(y_3=\infty\). Let \(x_1'=g'x_1\) and \(x_2'=g'x_2\). Then the matrix
  \(g''=\begin{pmatrix}
    1 &  -x_1'\\
    0 & 1
  \end{pmatrix}\in N\), fixes \(\infty=g' x_3\) and sends \(x'_1\) to
  \(y_1=0\). Let \(x''_2=g''x_2'\), from positive orientation of
  triples we have \(0<x''_2<\infty\). Next, the matrix
  \(g'''=\begin{pmatrix}
    a^{-1} &  0\\
    0 & a
  \end{pmatrix}\in A\) with \(a=\sqrt{x''_2}\)  sends \(x''_2\) to \(1\)
  and fixes both
  \(\infty=g''g'x_3\) and \(0=g''g' x_1\). Thus, \(g=g'''g''g'\) makes
  the required transformation \((x_1,x_2,x_3)\mapsto(0,1,\infty)\).
\end{proof}
\begin{corollary}
  Let \(\{x_1,x_2,x_3\}\) and \(\{y_1,y_2,y_3\}\) be two triples with the
  opposite orientations. Then, there is a
  unique M\"obius map \(\phi\in\SL\) with \(\phi\circ \cycle{}{0\infty}(x_j)=y_j\) for
  \(j=1\), \(2\), \(3\).
\end{corollary}
We will denote by \(\phi_{XY}\) the unique map from
Lem.~\ref{le:three-transitive} defined by triples
\(X=\{x_1,x_2,x_3\}\) and \(Y=\{y_1,y_2,y_3\}\).  

Although we are not going to use it in this paper, we note that the
following important result~\cite{Lang85}*{\S~III.1} is an immediate
consequence of our \emph{proof} of Lem.~\ref{le:three-transitive}.
\begin{corollary}[Iwasawa decomposition]
  \label{co:Iwasawa-decomp}
  Any element of \(g\in\SL\) is a product \(g=g_A g_N g_K\),  where
  \(g_A\),  \(g_N\) and  \(g_K\) belong to subgroups \(A\), \(N\),
  \(K\) respectively and those factors are uniquely defined.
\end{corollary}
In particular, we note that it is not a coincidence that the subgroups
appear in the Iwasawa decomposition \(\SL=ANK\) in order of decreasing
number of their fixed points.

\subsection{Triples of intervals}
\label{sec:triples-intervals}
We change our point of view and instead of two ordered triples of points
consider three ordered pairs, that is---three intervals or,
equivalently, three cycles. This is done in line with our extension of
Lie--M\"obius geometry.

For such triples of intervals we will need the following definition.
\begin{definition}
  \label{de:aligned-triples-intervals}
  We say that a triple of intervals \(\{[x_1,y_1], [x_2,y_2],
  [x_3,y_3]\}\)  is \emph{aligned} if the triples \(X=\{x_1,x_2,x_3\}\)
  and \(Y=\{y_1,y_2,y_3\}\)  of their endpoints have the same orientation. 
\end{definition}
Aligned triples determine certain one-parameter subgroups of M\"obius
transformations as follows:
\begin{proposition}
  \label{pr:triple-subgroup}
  Let \(\{[x_1,y_1], [x_2,y_2], [x_3,y_3]\}\) be an aligned triple of
  intervals.
  \begin{enumerate}
  \item If \(\phi_{XY}\) has at most one fixed point, then there is a
    unique (up to a parametrisation) 
    one-parameter semigroup of
    M\"obius map \(\psi(t)\subset\SL\), which maps \([x_1,y_1]\) to
    \([x_2,y_2]\) and \([x_3,y_3]\):
    \begin{displaymath}
      \psi(t_j)(x_1)=x_j,\quad \psi(t_j)(y_1)=y_j,\qquad \text{ for
        some } t_j\in\Space{R}{}  \text{ and } j=2,3.
    \end{displaymath}
  \item Let \(\phi_{XY}\) have two fixed points \(x<y\) and \(\cycle{}{xy}\)
    be the orientation inverting M\"obius transformation with the
    matrix~\eqref{eq:matrix-for-cycle-defn}. For \(j=1\), 
    \(2\), \(3\), we define:
    \begin{align*}
      x'_j&=x_j,& y'_j&=y_j,& x''_j&=\cycle{}{xy}x_j,&
      y''_j&=\cycle{}{xy}y_j&\text{ if } x<x_j<y,\\
      x'_j&=\cycle{}{xy}x_j,& y'_j&=\cycle{}{xy}y_j,& x''_j&=x_j,& y''_j&=y_j,&
      \text{ otherwise}.
    \end{align*}
    Then, there is a 
    one-parameter
    semigroup of M\"obius map \(\psi(t)\subset\SL\), and \(t_2\),
    \(t_3\in\Space{R}{}\) such that:
    \begin{displaymath}
      \psi(t_j)(x'_1)=x'_j,\quad \psi(t_j)(x''_1)=x''_j,\quad
      \psi(t_j)(y'_1)=y'_j,\quad  \psi(t_j)(y''_1)=y''_j\notingiq
    \end{displaymath}
    where \(j=2,3\).
  \end{enumerate}
\end{proposition}
\begin{proof}
  Consider the one-parameter subgroup of
  \(\psi(t)\subset \SL\) such that \(\phi_{XY}=\psi(1)\). Note, that
  \(\psi(t)\) and \(\phi_{XY}\) have the same fixed points (if any) and no
  point \(x_j\) is fixed since \(x_j\neq y_j\). If the number of fixed
  points is less than \(2\), then \(\psi(t)x_1\), \(t\in\Space{R}{}\)
  produces the entire real line except a possible single fixed
  point. Therefore, there are \(t_2\) and \(t_3\) such that
  \(\psi(t_2)x_1=x_2\) and \(\psi(t_3)x_1=x_3\). Since \(\psi(t)\) and
  \(\phi_{XY}\) commute for all \(t\) we also have:
  \begin{displaymath}
    \psi(t_j)y_1=\psi(t_j)\phi_{XY} x_1=\phi_{XY}\psi(t_j) x_1=\phi_{XY}
    x_j=y_j,\quad \text{ for } j=2,3.
  \end{displaymath}

  If there are two fixed points \(x<y\), then the open interval
  \((x,y)\) is an orbit for the subgroup \(\psi(t)\). Since all
  \(x'_1\), \(x'_2\) and \(x'_3\) belong to this orbit and \(\cycle{}{xy}\)
  commutes with \(\phi_{XY}\) we may repeat the above reasoning for
  the dashed intervals \([x'_j,y'_j]\). Finally, \(x''_j=\cycle{}{xy}x'_j\)
  and \(y''_j=\cycle{}{xy}y'_j\), where \(\cycle{}{xy}\) commutes with \(\phi\)
  and \(\psi(t_j)\), \(j=2\), \(3\).  Uniqueness of the subgroup
  follows from Lemma~\ref{le:uniqueness-orbit}.
\end{proof}
The group \(\SL\) acts transitively on collection of all cycles
\(\cycle{}{xy}\), thus this is a \(\SL\)-homogeneous space. It is easy to
see that the fix-group of the cycle \(\cycle{}{-1,1}\) is
\(\Aprime\)~\eqref{eq:a1n1-subgroup}. Thus the homogeneous space of
cycles is isomorphic to \(\SL/\Aprime\). 
\begin{lemma}
  \label{le:uniqueness-orbit}
  Let \(H\) be a one-parameter continuous subgroup of \(\SL\) and
  \(X=\SL/H\) be the corresponding homogeneous space. If two orbits
  of one-parameter continuous subgroups on \(X\) have at least three
  common points then these orbits coincide.
\end{lemma}
\begin{proof}
  Since \(H\) is conjugated either to \(\Aprime\), \(N'\) or \(K\),
  the homogeneous space \(X=\SL/H\) is isomorphic to the upper
  half-plane in double, dual or complex
  numbers~\cite{Kisil12a}*{\S~3.3.4}.  Orbits of one-parameter
  continuous subgroups in \(X\) are conic sections, which are circles,
  parabolas (with vertical axis) or equilateral hyperbolas (with
  vertical axis) for the respective type of geometry. Any two
  different orbits of the same type intersect at most at two points,
  since an analytic solution reduces to a quadratic equation. 
\end{proof}
Alternatively, we can reformulate Prop.~\ref{pr:triple-subgroup} as
follows: three different cycles \(\cycle{}{x_1y_1}\), \(\cycle{}{x_2y_2}\),
\(\cycle{}{x_3y_3}\) define a one-parameter subgroup, which generate either one
orbit or two related orbits passing the three cycles.

We have seen that the number of fixed points is the key
characteristics for the map \(\phi_{XY}\). The next result gives an
explicit expression for it.
\begin{proposition} 
\label{pr:type-of-subgroup}
  The map \(\phi_{XY}\) has zero, one or two fixed points
  if the expression
  \begin{equation}
    \label{eq:discriminant}
    \det
    \begin{pmatrix}
      1&x_1y_1&y_1-x_1\\
      1&x_2y_2&y_2-x_2\\
      1&x_3y_3&y_3-x_3
    \end{pmatrix}^2-4\det\begin{pmatrix}
      x_1&1&y_1\\
      x_2&1&y_2\\
      x_3&1&y_3
    \end{pmatrix}
    \cdot\det
    \begin{pmatrix}
      x_1&-x_1y_1&y_1\\
      x_2&-x_2y_2&y_2\\
      x_3&-x_3y_3&y_3
    \end{pmatrix}
  \end{equation}
  is negative, zero or positive respectively.
\end{proposition}
\begin{proof}
  If a M\"obius transformation \(\begin{pmatrix}
    a&b\\c&d
  \end{pmatrix}\) maps \(x_1\mapsto y_1\),  \(x_2\mapsto y_2\),
  \(x_3\mapsto y_3\) and  \(s\mapsto s\), then we have a homogeneous linear
  system, cf.~\cite{Beardon05a}*{Ex.~13.2.4}: 
  \begin{equation}
    \label{eq:three-pairs-fixed-point}
    \begin{pmatrix}
    x_1&1&-x_1y_1&-y_1\\
    x_2&1&-x_2y_2&-y_2\\
    x_3&1&-x_3y_3&-y_3\\
    s& 1 &-s^2&-s
  \end{pmatrix}
  \begin{pmatrix}
    a\\b\\c\\d
  \end{pmatrix}=
  \begin{pmatrix}
    0\\0\\0\\0
  \end{pmatrix}.
  \end{equation}
  A non-zero solution exists if the determinant of the \(4\times 4\)
  matrix is zero. Expanding it over the last row and 
  rearranging terms we obtain the quadratic equation for the fixed point
  \(s\):
  \begin{displaymath}
    s^2 \det
    \begin{pmatrix}
      x_1&1&y_1\\
      x_2&1&y_2\\
      x_3&1&y_3
    \end{pmatrix}
    +s\det
    \begin{pmatrix}
      1&x_1y_1&y_1-x_1\\
      1&x_2y_2&y_2-x_2\\
      1&x_3y_3&y_3-x_3
    \end{pmatrix}
    +\det
    \begin{pmatrix}
      x_1&-x_1y_1&y_1\\
      x_2&-x_2y_2&y_2\\
      x_3&-x_3y_3&y_3
    \end{pmatrix}
    =0.
  \end{displaymath}
  The value~\eqref{eq:discriminant} is the discriminant of this equation.
\end{proof}
\begin{remark}
  It is interesting to note, that the relation \(ax+b-cxy-dy=0\) used
  in~\eqref{eq:three-pairs-fixed-point} 
  can be stated as e-orthogonality of the cycle \(
  \begin{pmatrix}
    a&b\\c&d
  \end{pmatrix}\) and the isotropic bilinear form \(
  \begin{pmatrix}
    x&-xy\\1&-y
  \end{pmatrix}\).
\end{remark}
 If \(y=g_0\cdot x\) for some \(g_0\in
H_\tau \), then for any \(g\in H_\tau\) we also have
\(y_g=g_0\cdot x_g\), where \(x_g=g\cdot x_g\) and \(y_g=g\cdot
y_g\). Thus, we demonstrated the first part of the following result.
\begin{lemma}
  Let \(\tau =1\), \(0\) or \(-1\) and a real constant \(t\neq 0\) be such that
    \(1-\tau t^2>0\).
  \begin{enumerate} 
  \item The collections of intervals:
    \begin{equation}
      \label{eq:tau-invariant-cycles}
      I_{\tau,t} =\left\{[x,\textstyle \frac{x+\tau t}{tx+1}]\such x\in \Space{R}{} \right\}
    \end{equation}
    is preserved by the actions of subgroup \(H_\tau\). Any three
    different intervals from \(I_{\tau,t}\) define the subgroup
    \(H_\tau\) in the sense of Prop.~\ref{pr:triple-subgroup}.
  \item 
    All \(H_\tau\)-invariant bilinear forms compose the family \(P_{\tau,t}=\left\{
      \begin{pmatrix}
        a&\tau b\\
        b& a
      \end{pmatrix}
    \right\}\).
  \end{enumerate}
\end{lemma}
The family \(P_{\tau,t}\) consists of the eigenvectors of the
\(4\times 4\) matrix from~\eqref{eq:SL2-act-cycle-linear} with
suitably substituted entries. There is (up to a factor) exactly one
\(\tau\)-isotropic form in \(P_{\tau,t}\), namely \(\begin{pmatrix}
  1&\tau \\
  1& 1
\end{pmatrix}\). We denote this form by \(\alli\). It corresponds to
the point \((0,1)\in\Space{R}{2}\) as discussed after
Defn.~\ref{de:product-cycles}. We may say that the subgroup
\(H_\tau\) fixes the point \(\alli\), this will play an important
r\^ole below. 

\subsection{Geometrisation of cycles}
\label{sec:geom-cycl}

We return to the geometric version of the Poincar\'e extension
considered in Sec.~\ref{sec:geom-constr} in terms of cycles. 
Cycles of the form \(\begin{pmatrix}
  x&-x^2\\1&-x
\end{pmatrix}\) are \(\tau\)-isotropic for any \(\tau\) and are
parametrised by the point \(x\) of the real line. For a fixed
\(\tau\), the collection of all \(\tau\)-isotropic cycles is a larger
set containing the image of the real line from the previous
sentence. Geometrisation of this embedding is described in the
following result.
\begin{lemma}
  \begin{enumerate}
    \item The transformation \(x\mapsto \frac{x+\tau t}{tx+1}\) from the
    subgroup \(H_\tau\), which maps \(x\mapsto y\), corresponds to
    the value \(t=\frac{x-y}{xy-\tau}\).  
  \item The unique (up to a factor) bilinear form \(Q\) orthogonal
    to \(\cycle{}{xx}\), \(\cycle{}{yy}\) and \(\alli\) is
    \begin{displaymath}
      Q=\begin{pmatrix}
        \frac{1}{2}(x+y+xy-\tau)&-xy\\
        1&\frac{1}{2}(-x-y+xy-\tau)
      \end{pmatrix}.
    \end{displaymath}
  \item The defined above \(t\) and  \(Q\) are connected by the identity:
    \begin{equation}
      \label{eq:cosine-real-axis}
      \frac{\scalar[\tau]{Q}{\Space{R}{}}}{\sqrt{\modulus{\scalar[\tau]{Q}{Q}}}}
      =\frac{\tau}{\sqrt{\modulus{t^2-\tau}}}.
    \end{equation}
    Here, the real line is represented by the bilinear form \(\Space{R}{}=
    \begin{pmatrix}
     2^{-1/2}&0\\0&2^{-1/2} 
    \end{pmatrix}
    \) normalised such
    that \(\scalar[\tau]{\Space{R}{}}{\Space{R}{}}=\pm1\).
  \item For a cycle \(Q=
  \begin{pmatrix}
    l+n&-m\\k&-l+n
  \end{pmatrix}
  \), the value
  \(\frac{\scalar[e]{Q}{\Space{R}{}}}{\sqrt{\modulus{\scalar[e]{Q}{Q}}}}
  =-\frac{n}{\sqrt{\modulus{l^2+n^2-km}}}\)
  is equal to the cosine of the angle between the curve \(k(u^2+\tau
  v^2)-2lu-2nv+m=0\)~\eqref{eq:quadratics-equation} and the real line,
  cf.~\cite{Kisil12a}*{Ex.~5.23}.
  \end{enumerate}
\end{lemma}
\begin{proof}
  The first statement is verified by a short calculation. A form \(Q=
  \begin{pmatrix}
    l+n&-m\\k&-l+n
  \end{pmatrix}
  \) in the
  second statement may be calculated from the homogeneous system:
  \begin{displaymath}
    \begin{pmatrix}
      0&-2x&x^2&1\\
      0&-2y&y^2&1\\
      -2&0&-\tau&1
    \end{pmatrix}
    \begin{pmatrix}
      n\\l\\k\\m
    \end{pmatrix}=
    \begin{pmatrix}
      0\\0\\0
    \end{pmatrix}\notingiq
  \end{displaymath}
  which has the rank \(3\) if \(x\neq y\). The third statement
  can be checked by a calculation as well. Finally, the last item is a
  particular case of the more general statement as indicated. Yet, we
  can derive it here from the implicit derivative
  \(\frac{dv}{du}=\frac{ku-l}{n}\) of the function \(k(u^2+\tau
  v^2)-2lu-2nv+m=0\)~\eqref{eq:quadratics-equation} at the point
  \((u,0)\). Note that this value is independent from \(\tau\). Since
  this is the tangent of the intersection angle with the real line,
  the square of the cosine of this angle is:
  \begin{displaymath}
    \frac{1}{1+(\frac{dv}{du})^2}=\frac{n^2}{l^2+k^2 u^2+n^2-2 k u
      l}=\frac{n^2}{l^2+n^2-km}=
    \frac{\scalar{Q}{\Space{R}{}}^2}{\scalar[e]{Q}{Q}}\notingiq
  \end{displaymath}
  if \(ku^2-2  u l+m=0\).
\end{proof}
Also, we note that, the independence of the left-hand side
of~\eqref{eq:cosine-real-axis} from  \(x\) can be shown
from basic principles. Indeed, for a fixed \(t\) the subgroup
\(H_\tau\) acts transitively on the family of triples \(\{x,
\frac{x+\tau t}{tx+1}, \alli\}\), thus \(H_\tau\) acts transitively on
all bilinear forms orthogonal to such triples. However, the left-hand
side of~\eqref{eq:cosine-real-axis} is \(\SL\)-invariant, thus may not
depend on \(x\). This simple reasoning cannot provide the exact
expression in the right-hand side of~\eqref{eq:cosine-real-axis},
which is essential for the geometric interpretation of the Poincar\'e
extension.

To restore a cycle from its intersection points with the real line we
need also to know its cycle product with the real line. If this
product is non-zero then the sign of the parameter \(n\) is
additionally required.  At the cycles' language, a common point of
cycles \(\cycle{}{}\) and \(\cycle[\tilde]{}{}\) is encoded by a cycle
\(\hat{C}\) such that:
\begin{equation}
  \label{eq:cycles-passing-point}
  \scalar[e]{\cycle[\hat]{}{}}{\cycle{}{}}=
  \scalar[e]{\cycle[\hat]{}{}}{\cycle[\tilde]{}{}}=
  \scalar[\tau]{\cycle[\hat]{}{}}{\cycle[\hat]{}{}}=0.
\end{equation}
For a given value of \(\tau\), this produces two linear and one
quadratic equation for parameters of \(\hat{C}\). Thus, a pair of
cycles may not have a common point or have up to two such
points. Furthermore, M\"obius-invariance of the above
conditions~\eqref{eq:cosine-real-axis}
and~\eqref{eq:cycles-passing-point} supports the geometrical
construction of Poincar\'e extension, cf. Lem.~\ref{le:pencil-common-point}:
\begin{lemma}
  \label{le:poincare-geom-general}
  Let a family consist of cycles, which are \(e\)-orthogonal to a
  given \(\tau\)-isotropic cycle \(\cycle[\hat]{}{}\) and have the
  fixed value of the fraction in the left-hand side
  of~\eqref{eq:cosine-real-axis}. Then, for a given M\"obius
  transformation \(g\) and any cycle \(\cycle{}{}\) from the
  family, \(g\cycle{}{}\) is \(e\)-orthogonal to the
  \(\tau\)-isotropic cycle \(g\cycle[\hat]{}{}\) and has the same
  fixed value of the fraction in the left-hand side
  of~\eqref{eq:cosine-real-axis} as \(\cycle{}{}\).
\end{lemma}

Summarising the geometrical construction, the Poincar\'e extension based on
two intervals and the additional data produces two situations:
\begin{enumerate}
\item If the cycles \(C\) and \(C'\) are orthogonal to the real line, then
  a pair of overlapping cycles produces a point of the elliptic upper
  half-plane, a pair of disjoint cycles defines a point of the
  hyperbolic. However, there is no orthogonal cycles uniquely defining
  a parabolic extension.
\item If we admit cycles, which are not orthogonal to the real line,
  then the same pair of cycles may define any of the three different
  types (EPH) of extension.
\end{enumerate}
These peculiarities make the extension based on three intervals,
described above, a bit more preferable.

\subsection{Summary of the construction and generalisations}
\label{sec:concl-remarks-open}

Based on the consideration in
Sections~\ref{sec:cycles}--~\ref{sec:geom-cycl} we describe the
following steps to carry out the generalised extension procedure:
\begin{enumerate}
\item Points of the extended space are equivalence classes of aligned
  triples of cycles in \(P\Space{R}{1}\), see
  Defn.~\ref{de:aligned-triples-intervals}. The equivalence relation
  between triples will emerge at step~\ref{it:equivalence-triples}.
\item \label{it:one-param-subgroup} A triple \(T\) of different cycles
  defines the unique one-parameter continues subgroup \(S(t)\) of
  M\"obius transformations as defined in
  Prop.~\ref{pr:triple-subgroup}.
\item \label{it:equivalence-triples} Two triples of cycles are
  equivalent if and only if the subgroups defined in
  step~\ref{it:one-param-subgroup} coincide (up to a parametrisation).
\item The geometry of the extended space, defined by the equivalence
  class of a triple \(T\), is elliptic, parabolic or hyperbolic
  depending on the subgroup \(S(t)\) being similar \(S(t)=gH_\tau(t)
  g^{-1}\), \(g\in\SL\) (up to parametrisation) to
  \(H_\tau\)~\eqref{eq:a1n1k-universal} with \(\tau=-1\), \(0\) or
  \(1\) respectively. The value of \(\tau\) may be identified from the
  triple using Prop.~\ref{pr:type-of-subgroup}.
\item For the above \(\tau\) and \(g\in\SL\), the point of the
  extended space, defined by the the equivalence class of a triple
  \(T\), is represented by \(\tau\)-isotropic (see
  Defn.~\ref{de:product-cycles}(\ref{it:isotropic})) bilinear form
  \(g^{-1}
  \begin{pmatrix}
    1&\tau\\1&1
  \end{pmatrix}
  g\), which is \(S\)-invariant, see the end of
  Section~\ref{sec:triples-intervals}.
\end{enumerate}

Obviously, the above procedure is more complicated that the geometric
construction from Section~\ref{sec:geom-constr}. There are reasons for
this, as discussed in Section~\ref{sec:geom-cycl}: our procedure is
uniform and we are avoiding consideration of numerous subcases created
by an incompatible selection of parameters. Furthermore, our presentation
is aimed for generalisations to M\"obius transformations of moduli
over other rings. This can be considered as an analog of Cayley--Klein
geometries~\citelist{\cite{Yaglom79}*{Apps.~A--B} \cite{Pimenov65a}}.

It shall be rather straightforward to adopt the extension for
\(\Space{R}{n}\). M\"obius transformations in \(\Space{R}{n}\) are
naturally expressed as linear-fractional transformations in Clifford
algebras, cf.~\cite{Cnops02a} and
Sect.~\ref{sec:multi-dimens-vari}. There is a similar classification
of subgroups based on fixed points~\cites{Ahlfors85a,Zoll87} in
multidimensional case. The M\"obius invariant matrix presentation of
cycles \(\Space{R}{n}\) is already
known~\citelist{\cite{Cnops02a}*{(4.12)} \cite{FillmoreSpringer90a}
  \cite{Kisil14a}*{\S~5}}. Of course, it is necessary to enlarge the
number of defining cycles from \(3\) to, say, \(n+2\). This shall have
a connection with Cauchy--Kovalevskaya extension considered in
Clifford analysis~\cites{Ryan90a,Sommen85a}.  Naturally, a
consideration of other moduli and rings may require some more serious
adjustments in our scheme.

Our construction is based on the matrix presentations of cycles. This
techniques is effective in many different cases
\cites{Kisil12a,Kisil14a}. Thus, it is not surprising that such ideas
(with some technical variation) appeared independently on many
occasions~\citelist{\cite{Cnops02a}*{(4.12)}
  \cite{FillmoreSpringer90a} \cite{Schwerdtfeger79a}*{\S~1.1}
  \cite{Kirillov06}*{\S~4.2}}. The interesting feature of the present
work is the complete absence of any (hyper)complex numbers. It deemed
to be unavoidable \cite{Kisil12a}*{\S~3.3.4} to employ complex, dual
and double numbers to represent three different types of M\"obius
transformations extended from the real line to a plane. Also (hyper)complex
numbers were essential in~\cites{Kisil12a,Kisil06a} to define 
three possible types of cycle product~\eqref{eq:cycle-product-expl1},
and now we managed without them.

Apart from having real entries, our matrices for cycles share the
structure of matrices from~\citelist{\cite{Cnops02a}*{(4.12)}
  \cite{FillmoreSpringer90a} \cite{Kisil12a} \cite{Kisil06a}}. To
obtain another variant, one replaces the map
\(\map{i}\)~\eqref{eq:t-map-matrix} by
\begin{displaymath}
  \map{t}:\
  \begin{pmatrix}
    x_1&y_1\\
    x_2&y_2
  \end{pmatrix}
  \mapsto
  \begin{pmatrix}
    y_1&y_2\\
    x_1&x_2
  \end{pmatrix}\,.
\end{displaymath}
Then, we may define symmetric matrices in a manner similar
to~\eqref{eq:matrix-for-cycle-defn}: 
\begin{displaymath}
\cycle{t}{xy}=\frac{1}{2}M_{xy}\cdot\map{t}(M_{xy})=
  \begin{pmatrix}
    xy & \frac{x+y}{2}  \\
    \frac{x+y}{2} & 1
  \end{pmatrix}\,.
\end{displaymath}
This is the form of matrices for cycles similar
to~\citelist{\cite{Schwerdtfeger79a}*{\S~1.1}
  \cite{Kirillov06}*{\S~4.2}}. The property~\eqref{eq:covariance-Cxy}
with matrix similarity shall be replaced by the respective one with
matrix congruence: \(g\cdot \cycle{t}{xy} \cdot g^{t}=
\cycle{t}{x'y'}\). The rest of our construction may be adjusted
for these matrices accordingly.

\section{Example: Continued Fractions}
\label{sec:exampl-cont-fract}

Continued fractions remain an important and attractive topic of
current research \citelist{\cite{Khrushchev08a}
  \cite{BorweinPoortenShallitZudilin14a} \cite{Karpenkov2013a}
  \cite{Kirillov06}*{\S~E.3}}. A fruitful and geometrically appealing
method considers a continued fraction as an (infinite) product of
linear-fractional transformations from the M\"obius group, see
Sec.~\ref{sec:continued-fractions} below for a brief overview,
works~\citelist{\cite{PaydonWall42a} \cite{Schwerdtfeger45a}
  \cite{PiranianThron57a} \cite{Beardon04b}
  \cite{Schwerdtfeger79a}*{Ex.~10.2}} and in particular
\citelist{\cite{Beardon01a} \cite{SimonBasicComplex15}*{\S~7.5}}
contain further references and historical notes. Partial products of
linear-fractional maps form a sequence in the M\"obius group, the
corresponding sequence of transformations can be viewed as a discrete
dynamical system~\cite{Beardon01a,MageeOhWinter14a}. Many important
questions on continued fractions, e.g. their convergence, can be
stated in terms of asymptotic behaviour of the associated dynamical
system. Geometrical generalisations of continued fractions to many
dimensions were introduced recently as
well~\citelist{\cite{Beardon03a} \cite{Karpenkov2013a}}.

\subsection[Continued fractions and Moebius--Lie geometry]{Continued fractions and M\"obius--Lie geometry}
\label{sec:cont-fract-mobi}
A comprehensive consideration of the M\"obius group involves
cycles---the M\"obius invariant set of all circles and straight
lines. Furthermore, an efficient treatment cycles and M\"obius
transformations is realised through certain \(2\times 2\) matrices,
which we will review in Sec.~\ref{sec:mobi-transf-cycl}, see
also~\citelist{\cite{Schwerdtfeger79a} \cite{Cnops02a}*{\S~4.1}
  \cite{FillmoreSpringer90a} \cite{Kirillov06}*{\S~4.2}
  \cite{Kisil06a} \cite{Kisil12a}}. Linking the above relations we
present continued fractions within the extend M\"obius--Lie geometry:
\begin{claim}[Continued fractions and cycles]
  Properties of continued fractions may be illustrated and
  demonstrated using related cycles, in particular, in the form of
  respective \(2\times 2\) matrices. 
\end{claim}
One may expect that such an observation has been made a while ago,
e.g. in the book~\cite{Schwerdtfeger79a}, where both topics were
studied. However, this seems did not happen for some reasons. It is
only the recent paper~\cite{BeardonShort14a}, which pioneered a
connection between continued fractions and cycles. We hope that the
explicit statement of the claim will stimulate its further fruitful
implementations. In particular, Sec.~\ref{sec:cont-fract-cycl} reveals
all three possible cycle arrangements similar to one used
in~\cite{BeardonShort14a}. Secs.~\ref{sec:multi-dimens-vari}--\ref{sec:mult-cont-fract}
shows that relations between continued fractions and cycles can be
used in the multidimensional case as well.

As an illustration, we draw on Fig.~\ref{fig:cont-fract-e} chains of
tangent horocycles (circles tangent to the real line,
see~\cite{BeardonShort14a} and Sec.~\ref{sec:cont-fract-cycl}) for two
classical simple continued fractions:
\begin{displaymath}
  e=2+
  \cfrac{1}{1+\cfrac{1}{2+\cfrac{1}{1+\cfrac{1}{1+\ldots}}}}\,
  ,\qquad
  \pi=3+
  \cfrac{1}{7+\cfrac{1}{15+\cfrac{1}{1+\cfrac{1}{292+\ldots}}}}\, .
\end{displaymath}
\begin{figure}[htbp]
  \centering
  \includegraphics[scale=.6]{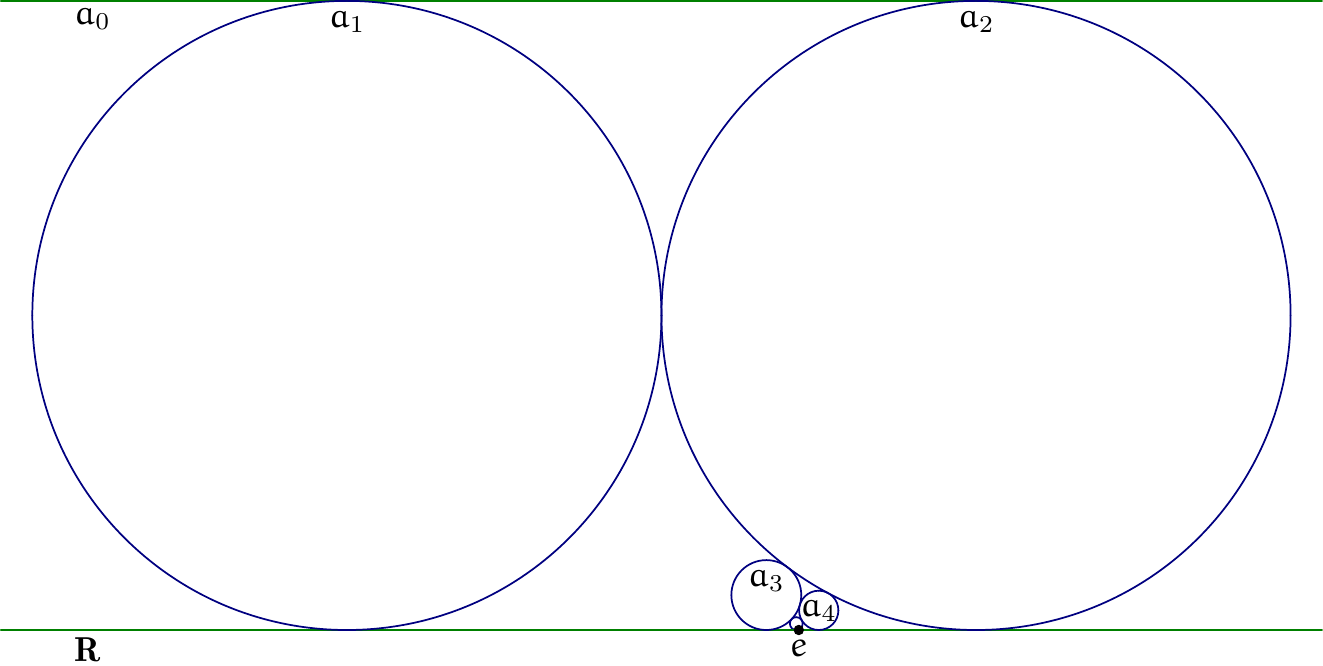}\hfill
  \includegraphics[scale=.6]{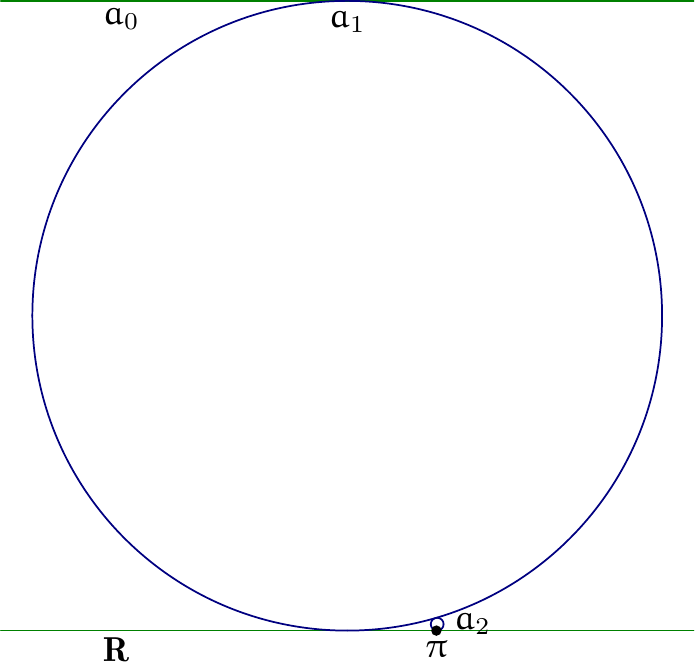}
  \caption{Continued fractions for $e$ and $\pi$ visualised. The
    convergence rate for $\pi$ is pictorially faster.}
  \label{fig:cont-fract-e}
\end{figure}
One can immediately see, that the convergence for \(\pi\) is much
faster: already the third horocycle is too small to be visible even if
it is drawn. This is related to larger coefficients in the continued
fraction for \(\pi\). 

Paper~\cite{BeardonShort14a} also presents examples of proofs based on
chains of horocycles.  This intriguing construction was introduced
in~\cite{BeardonShort14a} \emph{ad hoc}. Guided by the above claim we
reveal sources of this and other similar arrangements of horocycles.
Also, we can produce multi-dimensional variant of the framework.

\subsection{Preliminaries on Continued Fractions}
\label{sec:continued-fractions}

We use the following compact notation for a continued fraction:
\begin{equation}
\label{eq:cont-frac}
K(a_n|b_n)=
\cfrac{a_1}{b_1+\cfrac{a_2}{b_2+\cfrac{a_3}{b_3+\ldots}}}=
\frac{a_1}{b_1}\genfrac{}{}{0pt}{}{}{+}\frac{a_2}{b_2}\genfrac{}{}{0pt}{}{}{+}
\frac{a_3}{b_3}\genfrac{}{}{0pt}{}{}{+\ldots}.
\end{equation}
Without loss of generality we can assume \(a_j\neq 0\) for all
\(j\). 
The important particular case of \emph{simple} continued fractions,
with \(a_n=1\) for all \(n\), is denoted by \(K(b_n)=K(1|b_n)\). Any
continued fraction can be transformed to an equivalent simple one.

It is easy to see, that continued fractions are
related to the following linear-fractional
(M\"obius) transformation, cf.~\cites{PaydonWall42a,Schwerdtfeger45a,PiranianThron57a,Beardon04b}:
\begin{equation}
  \label{eq:cf-moebius-maps}
  S_n=s_1\circ s_2\circ \ldots\circ s_n, \qquad \text{where} \quad
  s_{j}(z)= \frac{a_j}{b_j+z}.
\end{equation}
These M\"obius transformation are considered as bijective maps of the
Riemann sphere \(\dot{\Space{C}{}}=\Space{C}{}\cup\{\infty\}\) onto itself.
If we associate  the matrix \(
\begin{pmatrix}
  a&b\\ c&d 
\end{pmatrix}\) to a liner-fractional transformation \(z\mapsto
\frac{a z+b}{c z+d}\), then the composition of two
such transformations corresponds to multiplication of the
respective matrices. Thus, relation~\eqref{eq:cf-moebius-maps}
has the matrix form:
\begin{equation}
  \label{eq:cf-part-frac-matrix}
  \begin{pmatrix}
    P_{n-1}&P_n\\
    Q_{n-1}&Q_n
  \end{pmatrix}
  =
  \begin{pmatrix}
    0&a_1\\
    1&b_1
  \end{pmatrix}
  \begin{pmatrix}
    0&a_2\\
    1&b_2
  \end{pmatrix}\cdots
  \begin{pmatrix}
    0&a_n\\
    1&b_n
  \end{pmatrix} .
\end{equation}
The last identity can be fold into the recursive formula: 
\begin{equation}
  \label{eq:cf-recurr-matrix}
  \begin{pmatrix}
    P_{n-1}&P_{n}\\    
    Q_{n-1}&Q_{n}
  \end{pmatrix}
  =
  \begin{pmatrix}
    P_{n-2}&P_{n-1}\\    
    Q_{n-2}&Q_{n-1}
  \end{pmatrix}
  \begin{pmatrix}
    0&a_n\\
    1&b_n
  \end{pmatrix}.
\end{equation}
This is equivalent to the main recurrence relation:
\begin{equation}
  \label{eq:cf-recurr}
  \begin{array}{c}
    P_n=b_nP_{n-1}+a_nP_{n-2}\\
    Q_n=b_nQ_{n-1}+a_nQ_{n-2}
  \end{array}, 
\quad n=1,2,3,\ldots,
\quad\text{with }
\begin{array}{cc}
{P_1}={a_1},&  {P_{0}}=0,\\
{Q_1}={b_1}, & Q_{0}=1.
\end{array}
\end{equation}

The meaning of entries \(P_n\) and \(Q_n\) from the
matrix~\eqref{eq:cf-part-frac-matrix} is revealed as follows. M\"obius
transformation~\eqref{eq:cf-moebius-maps}--\eqref{eq:cf-part-frac-matrix}
maps \(0\) and \(\infty\) to
\begin{equation}
  \label{eq:part-frac-moebius-map}
  \frac{P_n}{Q_n}=S_n(0), \qquad     \frac{P_{n-1}}{Q_{n-1}}=S_n(\infty).
\end{equation}
It is easy to see that \(S_n(0)\) is the \emph{partial quotient} of~\eqref{eq:cont-frac}:
\begin{equation}
  \label{eq:partial-quotient}
  \frac{P_n}{Q_n}= 
  \frac{a_1}{b_1}\genfrac{}{}{0pt}{}{}{+}\frac{a_2}{b_2}\genfrac{}{}{0pt}{}{}{+\ldots+}
\frac{a_n}{b_n}\,.
\end{equation}
Properties of the sequence of partial quotients
\(\left\{\frac{P_n}{Q_n}\right\}\) in terms of sequences \(\{a_n\}\)
and \(\{b_n\}\) are the core of the continued fraction
theory. Equation~\eqref{eq:part-frac-moebius-map} links partial
quotients with the M\"obius map~\eqref{eq:cf-moebius-maps}.  Circles
form an invariant family under M\"obius transformations, thus their
appearance for continued fractions is natural. Surprisingly, this
happened only recently in~\cite{BeardonShort14a}.

\subsection[Moebius Transformations and Cycles]{M\"obius Transformations and Cycles}
\label{sec:mobi-transf-cycl}

If \(M=\begin{pmatrix}
  a&b\\c&d
\end{pmatrix}\) is a matrix with real entries then for the purpose of
the associated M\"obius transformations \(M: z\mapsto
\frac{az+b}{cz+d}\) we may assume that \(\det M=\pm 1\). The
collection of all such matrices form a group.  M\"obius maps commute
with the complex conjugation \(z\mapsto \bar{z}\). If \(\det M>0\)
then both the upper and the lower half-planes are preserved. If \(\det
M <0\) then the two half-planes are swapped. Thus, we can treat \(M\)
as the map of equivalence classes \(z\sim\bar{z}\), which are labelled
by respective points of the closed upper half-plane. Under this
identification we consider any map produced by \(M\) with \(\det M
=\pm 1\) as the map of the closed upper-half plane to itself.

The characteristic property of M\"obius maps is that circles and lines
are transformed to circles and lines. We use the word \emph{cycles}
for elements of this M\"obius-invariant
family~\cites{Yaglom79,Kisil12a,Kisil06a}. We abbreviate a cycle given
by the equation
\begin{equation}
  \label{eq:cycle-defn}
  k(u^2+v^2)-2lv-2nu+m=0
\end{equation}
to the point \((k,l,n,m)\) of the three dimensional projective space
\(P\Space{R}{3}\). The equivalence relation \(z\sim \bar{z}\) is
lifted  to the equivalence relation 
\begin{equation}
  \label{eq:cycle-equiv}
  (k,l,n,m)\sim(k,l,-n,m)
\end{equation}
in the space of cycles, which again is compatible with the M\"obius
transformations acting on cycles.

The most efficient connection between cycles and M\"obius
transformations is realised through the construction, which may be
traced back to~\cite{Schwerdtfeger79a} and was subsequently
rediscovered by various authors~\citelist{\cite{Cnops02a}*{\S~4.1}
  \cite{FillmoreSpringer90a} \cite{Kirillov06}*{\S~4.2}}, see
also~\cites{Kisil06a,Kisil12a}. The construction associates a cycle
\((k,l,n,m)\) with the \(2\times 2\) matrix \(\cycle{}{}=\begin{pmatrix}
  l+\rmi  n&-m\\
  k&-l+\rmi n
\end{pmatrix}\), see discussion in~\cite{Kisil12a}*{\S~4.4} for a
justification. This identification is M\"obius covariant: the M\"obius
transformation defined by \(M=\begin{pmatrix} a&b\\c&d
\end{pmatrix}\) maps a cycle with matrix \(\cycle{}{}\) to the cycle with
matrix \(M\cycle{}{}M^{-1}\). Therefore, any M\"obius-invariant relation between cycles
can be expressed in terms of corresponding matrices. The central role
is played by the M\"obius-invariant inner pro\-duct~\cite{Kisil12a}*{\S~5.3}: 
\begin{equation}
  \label{eq:cycle-iner-pr}
  \scalar{\cycle{}{}}{\cycle[      \tilde]{}{}}=\Re \tr(\cycle{}{}\overline{\cycle[\tilde]{}{}})\notingiq
\end{equation}
which is a cousin of the product used in GNS construction of
\(C^*\)-algebras. Notably, the relation:
\begin{equation}
  \label{eq:first-ortho-gen}
  \scalar{\cycle{}{}}{\cycle[\tilde]{}{}}=0\quad \text{ or } \quad  km'+mk'-2nn'-2ll'=0\notingiq
\end{equation}
describes two cycles \(\cycle{}{}=(k,l,m,n)\) and
\(\cycle[\tilde]{}{}=(k',l',m',n')\)
orthogonal in Euclidean geometry.  Also, the inner
product~\eqref{eq:cycle-iner-pr} expresses the Descartes--Kirillov
condition~\citelist{\cite{Kirillov06}*{Lem.~4.1(c)}
  \cite{Kisil12a}*{Ex.~5.26}} of \(\cycle{}{}\) and
\(\cycle[\tilde]{}{}\) to be externally tangent:
\begin{equation}
  \label{eq:tangent-cycles}
  \scalar{\cycle{}{}+\cycle[\tilde]{}{}}{\cycle{}{}+\cycle[\tilde]{}{}} =0 
  \quad\text{or}\quad
  (l+l')^2+(n+n')^2 -(m+m')(k+k')=0\notingiq
\end{equation}
where the representing vectors \(\cycle{}{}=(k,l,n,m)\) and
\(\cycle[\tilde]{}{}=(k',l',m',n')\) from
\(P\Space{R}{3}\) need to be normalised by the conditions
\(\scalar{\cycle{}{}}{\cycle{}{}}=1\) and
\(\scalar{\cycle[\tilde]{}{}}{\cycle[\tilde]{}{}}=1\).

\subsection{Continued Fractions and Cycles}
\label{sec:cont-fract-cycl}

Now we are following~\cite{Kisil14a}.  Let
\(M=\begin{pmatrix} a&b\\c&d
\end{pmatrix}\) be a matrix with real entries and the determinant
\(\det M\) equal to \(\pm 1\), we denote this by \(\delta=\det M\). As
mentioned in the previous section, to calculate the image of a cycle
\(\cycle{}{}\) under M\"obius transformations \(M\) we can use matrix
similarity \(M\cycle{}{}M^{-1}\). If \(M=
\begin{pmatrix}
  P_{n-1}&P_n\\
  Q_{n-1}&Q_n 
\end{pmatrix}\) is the
matrix~\eqref{eq:cf-part-frac-matrix} associated to a continued
fraction and we are interested in the partial fractions
\(\frac{P_n}{Q_n}\), it is natural to ask: 
\begin{quote}
  \emph{Which cycles \(\cycle{}{}\) have transformations
    \(M\cycle{}{}M^{-1}\) depending on the first (or on the second)
    columns of \(M\) only?}
\end{quote}

It is a straightforward
calculation with matrices\footnote{This calculation can be done with
  the help of the tailored Computer Algebra System (CAS) as described
  in~\citelist{\cite{Kisil12a}*{App.~B}\cite{Kisil05b}}.} to check the
following statements: 

\begin{lemma}
  \label{le:first-col}
  The cycles \((0,0,1,m)\) (the horizontal lines \(v=m\)) are the only
  cycles, such that their
  images under the M\"obius transformation \(\begin{pmatrix}
    a&b\\c&d
  \end{pmatrix}\) are independent from the column \(
  \begin{pmatrix}
    b\\d
  \end{pmatrix}\). The image 
  associated to the column \(\begin{pmatrix} a\\c
  \end{pmatrix}\) is the horocycle \((c^2m,acm,\delta,a^2m)\), which
  touches the real line at \(\frac{a}{c}\) and has the radius
  \(\frac{1 }{mc^2}\).
\end{lemma}
In particular, for the matrix~\eqref{eq:cf-recurr-matrix} the
horocycle is touching the real line at the point
\(\frac{P_{n-1}}{Q_{n-1}}=S_n(\infty)\)~\eqref{eq:part-frac-moebius-map}. 
\begin{lemma}
  \label{le:second-col}
  The cycles \((k,0,1,0)\) (with the equation \(k(u^2+v^2)-2v=0\)) are
  the only cycles, such that their images under the M\"obius
  transformation \(\begin{pmatrix} a&b\\c&d
  \end{pmatrix}\) are independent from the column \(
  \begin{pmatrix}
    a\\c
  \end{pmatrix}\). The image 
  associated to the column  \(\begin{pmatrix} b\\d
  \end{pmatrix}\) is the horocycle
  \((d^2k,bdk,\delta,b^2k)\), which touches the real line at
  \(\frac{b}{d}\) and has the radius \(\frac{1 }{kd^2}\).
\end{lemma}
 In particular, for the matrix~\eqref{eq:cf-recurr-matrix} the horocycle is
touching the real line at the point
\(\frac{P_n}{Q_n}=S_n(0)\)~\eqref{eq:part-frac-moebius-map}. 
In view of these partial quotients the following cycles joining them
are of interest.
\begin{lemma}
  \label{le:blue}
  A cycle \((0,1,n,0)\) (any non-horizontal line passing
  \(0\)) is transformed by
  \eqref{eq:cf-moebius-maps}--\eqref{eq:cf-part-frac-matrix} to the
  cycle \((2Q_nQ_{n-1}, P_nQ_{n-1}+Q_nP_{n-1},\delta n,2P_nP_{n-1})\),
  which passes points \(\frac{P_n}{Q_n}=S_n(0)\) and
  \(\frac{P_{n-1}}{Q_{n-1}}=S(\infty)\) on the real line.
\end{lemma}
The above three families contain cycles with specific relations to
partial quotients through M\"obius transformations. There is one
degree of freedom in each family: \(m\), \(k\) and \(n\),
respectively. We can use the parameters to create an ensemble of three
cycles (one from each family) with some M\"obius-invariant
interconnections. Three most natural arrangements are illustrated by
Fig.~\ref{fig:vari-arrang-three}.  The first row presents the initial
selection of cycles, the second row---their images after a M\"obius
transformation (colours are preserved). The arrangements are as
follows:
\begin{figure}[htbp]
  \centering
  \includegraphics{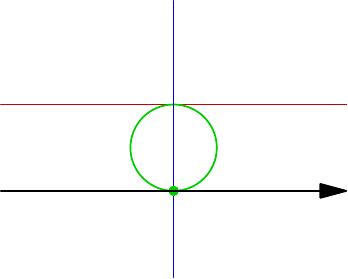}\hfil
  \includegraphics{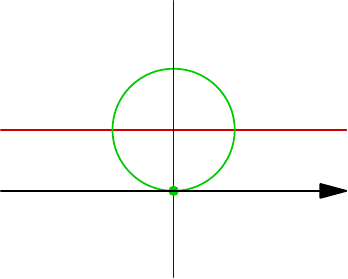}\hfil
  \includegraphics{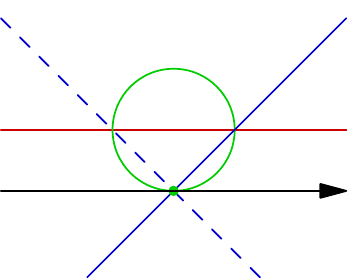}
\\
  \includegraphics{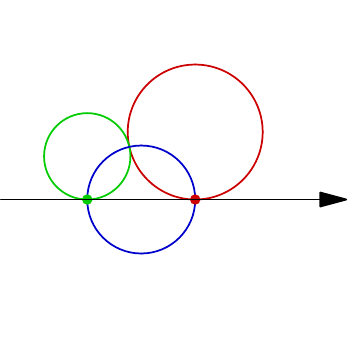}\hfil
  \includegraphics{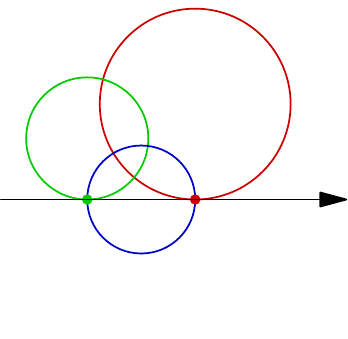}\hfil
  \includegraphics{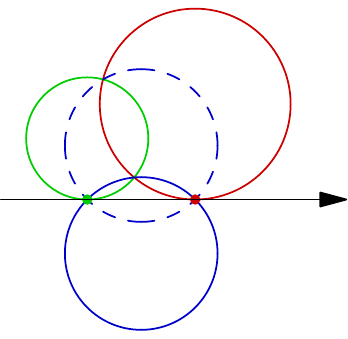}
  \caption[Various arrangements for three cycles]{Various arrangements
    for three cycles. The first row shows the initial position, the
    second row---after a M\"obius
    transformation (colours are preserved). \\
    The left column shows the arrangement used in the
    paper~\cite{BeardonShort14a}: two horocycles touching, the connecting
    cycle is passing their common point
    and is orthogonal to the real line.\\
    The central column presents two orthogonal horocycles and the
    connecting
    cycle orthogonal to them. \\
    The horocycles in the right column are again orthogonal but the
    connecting cycle passes one of their intersection points and makes
    \(45^\circ\) with the real axis.}
  \label{fig:vari-arrang-three}
\end{figure}

\begin{enumerate}
\item The left column shows the arrangement used in the
  paper~\cite{BeardonShort14a}: two horocycles are tangent, the third
  cycle, which we call \emph{connecting}, passes three points of
  pair-wise contact between horocycles and the real line. The
  connecting cycle is also orthogonal to horocycles and the real
  line. The arrangement corresponds to the following values \(m=2\),
  \(k=2\), \(n=0\). These parameters are uniquely defined by the above
  tangent and orthogonality conditions together with the requirement
  that the horocycles' radii agreeably depend from the consecutive
  partial quotients' denominators: \(\frac{1}{2Q_{n-1}^2}\) and
  \(\frac{1}{2Q_{n}^2}\) respectively. This follows from the explicit
  formulae of image cycles calculated in Lemmas~\ref{le:first-col}
  and~\ref{le:second-col}.

\item The central column of Fig.~\ref{fig:vari-arrang-three} presents
  two orthogonal horocycles and the connecting cycle orthogonal to
  them. The initial cycles have parameters \(m=\sqrt{2}\),
  \(k=\sqrt{2}\), \(n=0\). Again, these values follow from the
  geometric conditions and the alike dependence of radii from
  the partial quotients' denominators: \(\frac{\sqrt{2}}{2Q_{n-1}^2}\) and
  \(\frac{\sqrt{2}}{2Q_{n}^2}\).

\item Finally, the right column have the same two orthogonal
  horocycles, but the connecting cycle passes one of two horocycles'
  intersection points. Its mirror reflection in the real
  axis satisfying~\eqref{eq:cycle-equiv} (drawn in the dashed style) passes
  the second intersection point.  This corresponds to the values
  \(m=\sqrt{2}\), \(k=\sqrt{2}\), \(n=\pm 1\). The connecting cycle makes
  the angle \(45^\circ\) at the points of intersection with the real
  axis. It also has the radius
  \(\frac{\sqrt{2}}{2}\modulus{\frac{P_n}{Q_n}-\frac{P_{n-1}}{Q_{n-1}}}
  =\frac{\sqrt{2}}{2}\frac{1}{\modulus{Q_nQ_{n-1}}}\)---the geometric
  mean of radii of two other cycles. This again repeats the relation
  between cycles' radii in the first case.
\end{enumerate}
Three configurations have fixed ratio \(\sqrt{2}\) between respective
horocycles' radii. Thus, they are equally suitable for the proofs based on
the size of horocycles, e.g.~\cite{BeardonShort14a}*{Thm.~4.1}. 

On the other hand, there is a tiny computational advantage in the case
of orthogonal horocycles. Let we have the sequence \(p_j\) of partial
fractions \(p_j=\frac{P_j}{Q_j}\) and want to rebuild the
corresponding chain of horocycles. A horocycle with the point of
contact \(p_j\) has components \((1,p_j, n_j, p_j^2)\), thus only the
value of \(n_j\) need to be calculated at every step. If we use the
condition ``to be tangent to the previous horocycle'', then the quadratic
relation~\eqref{eq:tangent-cycles} shall be solved. Meanwhile, the
orthogonality relation~\eqref{eq:first-ortho-gen} is linear in \(n_j\).

\subsection{Multi-dimensional M\"obius maps and cycles}
\label{sec:multi-dimens-vari}

It is natural to look for multidimensional generalisations of
continued fractions.  A geometric approach based on M\"obius
transformation and Clifford algebras was proposed
in~\cite{Beardon03a}. 
Let us restrict the consideration from
Sect.~\ref{sec:cliff-algebr-mobi} to Clifford algebra \(\Cliff{n}\) of
the Euclidean space \(\Space{R}{n}\). It is the associative unital
algebra over \(\Space{R}{}\) generated by the elements
\(e_1\),\ldots,\(e_n\) satisfying the following relation:
\begin{displaymath}
  e_i e_j + e_je_i=-2\delta_{ij}\notingiq
\end{displaymath}
where \(\delta_{ij}\) is the Kronecker delta. An element of
\(\Cliff{n}\) having the form \(x=x_1e_1+\ldots+x_ne_n\) can be
associated with the vector \((x_1,\ldots,x_n)\in\Space{R}{n}\).  The
\emph{reversion} \(a\mapsto a^*\) in
\(\Cliff{n}\)~\cite{Cnops02a}*{(1.19(ii))} is defined on vectors by
\(x^*=x\) and extended to other elements by the relation
\((ab)^*=b^*a^*\). Similarly the \emph{conjugation} is defined on
vectors by \(\bar{x}=-x\) and the relation
\(\overline{ab}=\bar{b}\bar{a}\). We also use the notation
\(\modulus{a}^2=a\bar{a}\geq 0\) for any product \(a\) of vectors.
An important observation is that any non-zero vectors \(x\) has a
multiplicative inverse: \(x^{-1}=\frac{\bar{x}}{\modulus{x}^2}\).

By Ahlfors~\cite{Ahlfors86} (see
also~\citelist{\cite{Beardon03a}*{\S~5} \cite{Cnops02a}*{Thm.~4.10}})
a matrix \(M=
\begin{pmatrix}
  a&b\\c&d
\end{pmatrix}\) with Clifford entries defines a linear-fractional transformation of
\(\Space{R}{n}\) if the
following conditions are satisfied:
\begin{enumerate}
\item \(a\), \(b\), \(c\) and \(d\) are products of vectors in
  \(\Space{R}{n}\){\notingiqsemitocoma}
\item\label{it:ab-cd-ca-db-vectors} \(ab^*\), \(cd^*\), \(c^*a\) and \(d^*b\) are vectors in
  \(\Space{R}{n}\){\notingiqsemitocoma}
\item the pseudodeterminant \(\delta:=ad^*-bc^*\) is a non-zero real number.
\end{enumerate}
Clearly we can scale the matrix to have the pseudodeterminant
\(\delta=\pm 1\) without an effect on the related linear-fractional
transformation. Recall notations~\eqref{eq:matrix-bar-star}, cf.~\cite{Cnops02a}*{(4.7)}
\begin{displaymath}
  \bar{M}=
\begin{pmatrix}
  d^*&-b^*\\-c^*&a^*
\end{pmatrix}\qquad \text{ and } \qquad
M^*=\begin{pmatrix}
  \bar{d} &\bar{b}\\\bar{c}&\bar{a}
\end{pmatrix}.
\end{displaymath}
Then \(M\bar{M}=\delta I\) and \(\bar{M}=\kappa M^*\), where
\(\kappa=1\) or \(-1\) depending either \(d\) is a product of even or
odd number of vectors.

To adopt the discussion from Section~\ref{sec:mobi-transf-cycl} to
several dimensions we use vector rather than paravector
formalism, see~\cite{Cnops02a}*{(1.42)} for a discussion. 
Namely, we consider vectors \(x\in\Space{R}{n+1}\) as elements
\(x=x_1e_1+\ldots+x_ne_n+x_{n+1} e_{n+1}\) in \(\Cliff{n+1}\). 
Therefore we can extend the M\"obius
transformation defined by \(M=
\begin{pmatrix}
  a&b\\c&d
\end{pmatrix}\) with \(a,b,c,d\in\Cliff{n}\) to act on
\(\Space{R}{n+1}\). Again, such transformations commute with the
reflection \(R\) in the hyperplane \(x_{n+1}=0\):
\begin{displaymath}
  R:\quad x_1e_1+\ldots+x_ne_n+x_{n+1} e_{n+1}\quad
  \mapsto \quad x_1e_1+\ldots+x_ne_n-x_{n+1} e_{n+1}.
\end{displaymath}
Thus we can consider the M\"obius maps acting on the equivalence
classes \(x\sim R(x)\).

Spheres and hyperplanes in \(\Space{R}{n+1}\)---which we continue to
call cycles---can be associated to \(2\times 2\) matrices
by~\eqref{eq:spheres-Rn}, cf.~\citelist{\cite{FillmoreSpringer90a}
  \cite{Cnops02a}*{(4.12)}}:
\begin{displaymath}
  k\bar{x}x-l\bar{x}-x\bar{l}+m=0 \quad \leftrightarrow \quad 
  \cycle{}{}=
  \begin{pmatrix}
    l & m\\
    k & \bar{l}
  \end{pmatrix}\notingiq
\end{displaymath}
where \(k, m\in\Space{R}{}\) and \(l\in\Space{R}{n+1}\).  For brevity
we also encode a cycle by its coefficients \((k,l,m)\). 
The identification is also M\"obius-covariant in the sense that the
transformation associated with the Ahlfors matrix \(M\) sends a cycle
\(\cycle{}{}\) to the cycle \(M\cycle{}{}M^{*}\)~\cite{Cnops02a}*{(4.16)}.
The equivalence  \(x\sim R(x)\) is extended to spheres:
\begin{displaymath}
  \begin{pmatrix}
    l & m\\
    k & \bar{l}
  \end{pmatrix}\quad\sim   \quad
  \begin{pmatrix}
    R(l) & m\\
    k & R(\bar{l})
  \end{pmatrix}
\end{displaymath}
since it is preserved by the M\"obius transformations with
coefficients from \(\Cliff{n}\). 

Similarly to~\eqref{eq:cycle-iner-pr} we define the M\"obius-invariant
inner product of cycles by the identity
\(\scalar{\cycle{}{}}{\cycle[\tilde]{}{}}=\Re
\tr(\cycle{}{}\cycle[\tilde]{}{})\), where \(\Re\) denotes the scalar
part of a Clifford number.  The orthogonality condition
\(\scalar{\cycle{}{}}{\cycle[\tilde]{}{}}=0\) means that the
respective cycle are geometrically orthogonal in
\(\Space{R}{n+1}\). 

\subsection{Continued fractions from Clifford algebras and horocycles}
\label{sec:mult-cont-fract}

There is an association between the triangular matrices and the elementary
M\"obius maps of \(\Space{R}{n}\), cf.~\eqref{eq:cf-moebius-maps}:
\begin{equation}
  \label{eq:clifford-inversion}
  \begin{pmatrix}
    0&1\\1&b
  \end{pmatrix}: \
  x\ \mapsto\ (x+b)^{-1}\,,\ 
  \text{ where }
  x=x_1e_1+\ldots+x_ne_n, b=b_1e_1+\ldots b_ne_n.
\end{equation}
Similar to the real line case in
Section~\ref{sec:continued-fractions}, Beardon
proposed~\cite{Beardon03a} to consider the composition of a series of
such transformations as multidimensional continued fraction,
cf.~\eqref{eq:cf-moebius-maps}. It can be again represented as the the
product~\eqref{eq:cf-part-frac-matrix} of the respective \(2\times 2\)
matrices. Another construction of multidimensional continued fractions
based on horocycles was hinted in~\cite{BeardonShort14a}. We wish to
clarify the connection between them.  The bridge is provided by the
respective modifications of Lem.~\ref{le:first-col}--\ref{le:blue}.
\begin{lemma}
  \label{le:first-col-mul}
  The cycles \((0,e_{n+1},m)\) (the ``horizontal'' hyperplane
  \(x_{n+1}=m\)) are the only cycles, such that their images under the
  M\"obius transformation \(\begin{pmatrix} a&b\\c&d
  \end{pmatrix}\) are independent from the column \(
  \begin{pmatrix}
    b\\d
  \end{pmatrix}\). The image 
  associated to the column \(\begin{pmatrix} a\\c
  \end{pmatrix}\) is the horocycle \((-m\modulus{c}^2, -ma\bar{c}+\delta  e_{n+1},m\modulus{a}^2)\), which
  touches the hyperplane \(x_{n+1}=0\) at \(\frac{a\bar{c}}{\modulus{c}^2}\) and has the radius
  \(\frac{1 }{m\modulus{c}^2}\).
\end{lemma}
\begin{lemma}
  \label{le:second-col-mul}
  The cycles \((k,e_{n+1},0)\) (with the equation \(k(u^2+v^2)-2v=0\)) are
  the only cycles, such that their images under the M\"obius
  transformation \(\begin{pmatrix} a&b\\c&d
  \end{pmatrix}\) are independent from the column \(
  \begin{pmatrix}
    a\\c
  \end{pmatrix}\). The image 
  associated to the column  \(\begin{pmatrix} b\\d
  \end{pmatrix}\) is the horocycle
  \((k\modulus{d}^2,kb\bar{d}+\delta e_{n+1},-kb \bar{b})\), which touches the
  hyperplane \(x_{n+1}=0\) at \(\frac{b\bar{d}}{\modulus{d}^2}\) and has
  the radius \(\frac{1 }{k\modulus{d}^2}\).
\end{lemma}
The proof of the above lemmas are reduced to multiplications of
respective matrices with Clifford entries. 
\begin{lemma}
\label{le:blue-multi}
A cycle \(\cycle{}{}=(0,l,0)\), where \(l=x+re_{n+1}\) and \(0\neq
x\in\Space{R}{n}\), \(r\in\Space{R}{}\), that is any non-horizontal
hyperplane passing the origin, is transformed into
\(M\cycle{}{}M^*=(c x\bar{d}+d\bar{x}\bar{c}, a
x\bar{d}+b\bar{x}\bar{c}+\delta r e_{n+1}, a x\bar{b}+
b\bar{x}\bar{a})\).  This cycle passes points
\(\frac{a\bar{c}}{\modulus{c}^2}\) and \(\frac{b\bar{d}}{\modulus{d}^2}\).  

If \(x= \bar{c}d\), then the centre of
\begin{displaymath}
  M\cycle{}{}M^*=(2\modulus{c}^2\modulus{d}^2,
a\bar{c}\modulus{d}^2+b\bar{d}\modulus{c}^2, (a\bar{c})(d\bar{b})+(b\bar{d})(c\bar{a}))
\end{displaymath}
is
\(\frac{1}{2}(\frac{a\bar{c}}{\modulus{c}^2}+\frac{b\bar{d}}{\modulus{d}^2})+
\frac{\delta r}{2\modulus{c}^2\modulus{d}^2}e_{n+1}\), that is, the centre
belongs to the two-dimensional plane passing the points
\(\frac{a\bar{c}}{\modulus{c}^2}\) and \(\frac{b\bar{d}}{\modulus{d}^2}\)
and orthogonal to the hyperplane \(x_{n+1}=0\).
\end{lemma}
\begin{proof}
  We note that \(e_{n+1} x=-x e_{n+1}\) for all
  \(x\in\Space{R}{n}\). Thus, for a product of  vectors \(d\in\Cliff {n}\)
  we have \(e_{n+1} \bar{d}=d^* e_{n+1}\). Then  
  \begin{displaymath}
    c e_{n+1}\bar{d}+d\bar{e}_{n+1}\bar{c}=(cd^*-dc^*)e_{n+1}=(cd^*-(cd^*)^*)e_{n+1}=0\notingiq
  \end{displaymath}
  due to the Ahlfors
  condition~\ref{it:ab-cd-ca-db-vectors}. Similarly, \(a
  e_{n+1}\bar{b}+ b\bar{e}_{n+1}\bar{a}=0\) and \(a
  e_{n+1}\bar{d}+b\bar{e}_{n+1}\bar{c}=(ad^*-bc^*) e_{n+1}=\delta e_{n+1}\). 

  The image \(M\cycle{}{}M^*\) of the cycle \(\cycle{}{}=(0,l,0)\) is
  \((c l\bar{d}+d\bar{l}\bar{c}, a l\bar{d}+b\bar{l}\bar{c}, a
  l\bar{b}+ b\bar{l}{a})\).  From the above calculations for
  \(l=x+re_{n+1}\) it becomes \((c x\bar{d}+d\bar{x}\bar{c}, a
  x\bar{d}+b\bar{x}\bar{c}+\delta r e_{n+1}, a x\bar{b}+
  b\bar{x}\bar{a})\). The rest of statement is verified by the substitution.
\end{proof}

Thus, we have exactly the same freedom to choose representing
horocycles as in Section~\ref{sec:cont-fract-cycl}: make two
consecutive horocycles either tangent or orthogonal.  To visualise
this, we may use the two-dimensional plane \(V\) passing the points of
contacts of two consecutive horocycles and orthogonal to
\(x_{n+1}=0\). It is natural to choose the connecting cycle (drawn in
blue on Fig.~\ref{fig:vari-arrang-three}) with the centre belonging to
\(V\), this eliminates excessive degrees of freedom. The corresponding
parameters are described in the second part of
Lem.~\ref{le:blue-multi}.
Then, the intersection of horocycles with \(V\)
are the same as on Fig.~\ref{fig:vari-arrang-three}. 

Thus, the continued fraction with the partial quotients
\(\frac{P_n\bar{Q}_n}{\modulus{Q_n}^2}\in\Space{R}{n}\) can be
represented by the chain of tangent/orthogonal horocycles. The
observation made at the end of Section~\ref{sec:cont-fract-cycl} on
computational advantage of orthogonal horocycles remains valid in
multidimensional situation as well.

As a further alternative we may shift the focus from horocycles to the
connecting cycle (drawn in blue on
Fig.~\ref{fig:vari-arrang-three}). The part of the space
\(\Space{R}{n}\) encloses inside the connecting cycle is the image
under the corresponding M\"obius transformation of the half-space of
\(\Space{R}{n}\) cut by the hyperplane \((0,l,0)\) from
Lem.~\ref{le:blue-multi}. Assume a sequence of connecting cycles
\(\cycle{}{j}\) satisfies the following two conditions, e.g. in
Seidel--Stern-type theorem~\cite{BeardonShort14a}*{Thm~4.1}:
\begin{enumerate}
\item for any \(j\), the cycle \(\cycle{}{j}\) is enclosed within the cycle
  \(\cycle{}{j-1}\){\notingiqsemitocoma}
\item the sequence of radii of \(\cycle{}{j}\) tends to zero.
\end{enumerate}
Under the above assumption the sequence of partial fractions
converges. Furthermore, if we use the connecting cycles in the third
arrangement, that is generated by the cycle  \((0,x+e_{n+1},0)\),
where \(\norm{x}=1\), \(x\in\Space{R}{n}\), then the above second condition
can be replaced by following 
\begin{enumerate}
\item[(\(2'\))] the sequence of \(x_{n+1}^{(j)}\) of \((n+1)\)-th
  coordinates of the centres of the connecting cycles \(\cycle{}{j}\)
  tends to zero.
\end{enumerate}
Thus, the sequence of connecting cycles is a useful tool to describe
a continued fraction even without a relation to horocycles.

Summing up, we started from multidimensional continued fractions
defined by the composition of M\"obius transformations in
Clifford algebras and associated to it the respective chain of
horocycles. This establishes the equivalence of two approaches proposed
in~\cites{Beardon03a} and \cite{BeardonShort14a} respectively.

\section[Extension of Moebius--Lie Geometry by Ensembles of Interrelated Cycles]{Extension of M\"obius--Lie Geometry by Ensembles of Interrelated Cycles}
\label{sec:library-figure}

Previous consideration suggests a far-reaching generalisation
M\"obius--Lie geometry which we are presenting in this section.

\subsection{Figures as ensembles of cycles}
\label{sec:figures-as-families}

We start from some examples of ensembles of cycles, which conveniently
describe FLT-invariant families of objects.

\begin{example}
  \label{ex:ensamble-math}
  \begin{enumerate}
  \item \label{it:poincare-extension} The Poincar\'e extension of
    M\"obius transformations from the real line to the upper
    half-plane of complex numbers is described by a pair of distinct
    intersecting cycles \(\{\cycle{}{1}, \cycle{}{2}\}\), that is
    \(\scalar{\cycle{}{1}}{\cycle{}{2}}^2\leq
    \scalar{\cycle{}{1}}{\cycle{}{1}}
    \scalar{\cycle{}{2}}{\cycle{}{2}}\).

    Two such pairs
    \(\{\cycle{}{1}, \cycle{}{2}\}\) and \(\{\cycle[\tilde]{}{1},
    \cycle[\tilde]{}{2}\}\) are equivalent if any cycle
    \(\cycle[\tilde]{}{i}\) from the
    second pair is a linear combination of (i.e. belongs to the pencil
    spanned by) the cycles from the first pair.

    A modification from Sect.~\ref{sec:triples-intervals} with
    ensembles of three cycles describes an extension from the real
    line to the upper half-plane of complex, dual or double numbers.
    The construction can be generalised to arbitrary dimensions,
    cf. Sect.~\ref{sec:concl-remarks-open}.
\item\label{it:param-loxodromes}
  A parametrisation of loxodromes is provided by a triple of
  cycles \(\{\cycle{}{1}, \cycle{}{2}, \cycle{}{3}\}\) such
  that, cf.~\cite{KisilReid18a} and Fig.~\ref{fig:equiv-param-loxodr}: 
  \begin{enumerate}
  \item \(\cycle{}{1}\) is orthogonal to \(\cycle{}{2}\) and  \(\cycle{}{3}\){\notingiqsemitocoma}
  \item \(\scalar{\cycle{}{2}}{\cycle{}{3}}^2\geq
    \scalar{\cycle{}{2}}{\cycle{}{2}}
    \scalar{\cycle{}{3}}{\cycle{}{3}}\).
  \end{enumerate}
  Then, main invariant properties of M\"obius--Lie geometry,
  e.g. tangency of loxodromes, can be expressed in terms of this
  parametrisation~\cite{KisilReid18a}. 
\item A continued fraction is described by an infinite ensemble of
  cycles \((\cycle{}{k})\) such that, cf. Sect.~\ref{sec:cont-fract-cycl}:
  \begin{enumerate}
  \item All \(\cycle{}{k}\) are touching the real line (i.e. are
    \emph{horocycles}){\notingiqsemitocoma}
  \item \((\cycle{}{1})\) is a horizontal line passing through
    \((0,1)\){\notingiqsemitocoma}
  \item \(\cycle{}{k+1}\) is tangent to \(\cycle{}{k}\) for all \(k>1\).
  \end{enumerate}
  Several similar arrangements are described in
  Sect.~\ref{sec:cont-fract-cycl}. The key analytic properties of
  continued fractions---their convergence---can be linked to
  asymptotic behaviour of such an infinite
  ensemble~\cite{BeardonShort14a}.
\item A remarkable relation exists between discrete integrable systems
  and M\"obius geometry of finite configurations of
  cycles~\cites{BobenkoSchief18a,%
    KonopelchenkoSchief02a,KonopelchenkoSchief02b,%
    KonopelchenkoSchief05a,SchiefKonopelchenko09a}.  It comes from
  ``reciprocal force diagrams'' used in 19th-century statics, starting
  with J.C.~Max\-well. It is demonstrated in that the geometric
  compatibility of reciprocal figures corresponds to the algebraic
  compatibility of linear systems defining these configurations. On
  the other hand, the algebraic compatibility of linear systems lies
  in the basis of integrable systems. In
  particular~\cites{KonopelchenkoSchief02a,KonopelchenkoSchief02b},
  important integrability conditions encapsulate nothing but a
  fundamental theorem of ancient Greek geometry.
\item \label{it:wave-envelope}
  An important example of an infinite ensemble is provided by the
  representation of an arbitrary wave as the envelope of a continuous
  family of spherical waves. A finite subset of spheres can be used as
  an approximation to the infinite family. Then, discrete snapshots of
  time evolution of sphere wave packets represent a FLT-covariant
  ensemble of cycles~\cite{Bateman55a}. Further physical applications
  of FLT-invariant ensembles may be looked at~\cite{Kastrup08a}.
\end{enumerate}
\end{example}

One can easily note that the above parametrisations of some objects by
ensembles of cycles are not necessary unique. Naturally, two ensembles
parametrising  the same object are again connected by
FLT-invariant conditions. We presented only one example
here, cf.~\cite{KisilReid18a}.
\begin{figure}[htbp]
  \centering
    \animategraphics[controls=true,width=.9\textwidth,poster=first]{50}{_loxodromes}{1}{200}
    \caption[Equivalent parametrisation of a loxodrome]{Animated
      graphics of equivalent three-cycle parametrisations of a
      loxodrome. The green cycle is \(\cycle{}{1}\), two red circles are
      \(\cycle{}{2}\) and \(\cycle{}{3}\).}
  \label{fig:equiv-param-loxodr}
\end{figure}
\begin{example}
  Two non-degenerate triples \(\{\cycle{}{1},\cycle{}{2},\cycle{}{3}\}\) and
  \(\{\cycle[\tilde]{}{1},\cycle[\tilde]{}{2},\cycle[\tilde]{}{3}\}\) parameterise the same loxodrome as
  in Ex.~\ref{ex:ensamble-math}\ref{it:param-loxodromes} if and only if all the following
  conditions are satisfied:
  \begin{enumerate}
  \item \label{item:same-pencil}
    Pairs \(\{\cycle{}{2},\cycle{}{3}\}\) and \(\{\cycle[\tilde]{}{2},\cycle[\tilde]{}{3}\}\)  span the same
    hyperbolic pencil. That is cycles \(\cycle[\tilde]{}{2}\) and \(\cycle[\tilde]{}{3}\) are linear
    combinations of \(\cycle{}{2}\) and \(\cycle{}{3}\) and vise versa.
  \item \label{item:same-lambda}
    Pairs \(\{\cycle{}{2},\cycle{}{3}\}\) and \(\{\cycle[\tilde]{}{2},\cycle[\tilde]{}{3}\}\) have the same
    normalised cycle product~\eqref{eq:norm-cycle-prod}:
    \begin{equation}
      \label{eq:equal-lambdas}
      \nscalar {\cycle{}{2}}{\cycle{}{3}}=\nscalar {\cycle[\tilde]{}{2}}{\cycle[\tilde]{}{3}}.
    \end{equation}
  \item \label{item:ellipt-hyperb-ident}
    The elliptic-hyperbolic identity holds:
    \begin{equation}
      \label{eq:ellipt-hyperb-equat}
      \frac{\arccosh\nscalar {\cycle{}{j}}{\cycle[\tilde]{}{j}}}{\arccosh\nscalar{\cycle{}{2}}{\cycle{}{3}}}
      \equiv
      \frac{1}{2\pi}\arccos\nscalar {\cycle{}{1}}{\cycle[\tilde]{}{1}} \pmod{1}\,\notingiq
    \end{equation}
    where \(j\) is either \(2\) or \(3\).
  \end{enumerate}
  Various triples of cycles parametrising the same loxodrome are
  animated on Fig.~\ref{fig:equiv-param-loxodr}.
\end{example}
The respective equivalence relation for parametrisation of elliptic,
parabolic and hyperbolic Poincar\'e
extensions,
cf. Ex.~\ref{ex:ensamble-math}(\ref{it:poincare-extension}), by triple
of cycles can be deduced from Prop.~\ref{pr:triple-subgroup}.

\subsection[Extension of Moebius--Lie geometry and its implementation]{Extension of M\"obius--Lie geometry and its implementation through functional approach}
\label{sec:extension-mobius-lie}

The above examples suggest that one can
expand the subject and applicability of M\"obius--Lie geometry through
the following formal definition.
\begin{definition}
  \label{de:extended-Lie-Mobius}
  Let \(X\) be a set, \(R \subset X\times X\) be an oriented graph
  on \(X\) 
  and \(f\) be a function on \(R\) with values in FLT-invariant
  relations from \S~\ref{sec:conn-quadr-cycl}. Then
  \emph{\((R,f)\)-ensemble} is a collection of cycles
  \(\{C_j\}_{j\in X}\) such that
  \begin{displaymath}
    C_i \text { and } C_j \text{ are in the relation } f(i,j) \text {
      for all } (i,j)\in R.
  \end{displaymath}
  For a fixed FLT-invariant equivalence 
  relations \(\sim\) on the set \(\mathcal{E}\) of all \((R,f)\)-ensembles,
  \emph{the extended M\"obius--Lie geometry} studies properties of cosets
  \(\mathcal{E}/\sim\).
\end{definition}
This definition can be suitably modified for
\begin{enumerate}
\item ensembles with relations of more then two cycles, and/or
\item ensembles parametrised by continuous sets \(X\), cf. wave
  envelopes in Ex.~\ref{ex:ensamble-math}\ref{it:wave-envelope}.
\end{enumerate}

The above extension was developed along with the realisation
the library {{\bf{}figure}} within the \emph{functional programming}
framework. More specifically, an object from the {{\bf{}class} \ {\bf{}figure}} stores
defining relations, which link new cycles to the previously introduced
ones. This also may be treated as classical geometric
compass-and-straightedge constructions, where new lines or circles are
drawn through already existing elements. If requested, an explicit
evaluation of cycles\vtick\ parameters from this data may be
attempted.

To avoid ``chicken or the egg'' dilemma all cycles are stored in a
hierarchical structure of generations, numbered by integers. The basic
principles are:
\begin{enumerate}
\item Any explicitly defined cycle (i.e., a cycle which is not related to any
  previously known cycle) is placed into generation-0{\notingiqsemitocoma}
\item Any new cycle defined by relations to \emph{previous} cycles
  from generations \(k_1\), \(k_2\), \ldots, \(k_n\) is placed to the
  generation \(k\) calculated as:
  \begin{equation}
    \label{eq:generation-calculation}
    k=\max(k_1,k_2,\ldots,k_n)+1 .
  \end{equation}
  This rule does not forbid a cycle to have a relation to itself,
  e.g. isotropy (self-orthogonality) condition
  \(\scalar{\cycle{}{}}{\cycle{}{}}=0\), which specifies point-like
  cycles, cf. relation~\ref{it:point-zero-radius} in~\S~\ref{sec:conn-quadr-cycl}.  In fact, this is the only
  allowed type of relations to cycles in the same (not even speaking
  about younger) generations.
\end{enumerate}
There are the following alterations of the above rules:
\begin{enumerate}
\item From the beginning, every figure has two pre-defined cycles: the
  real line (hyperplane) \(\cycle{}{\Space{R}{}}\), and the zero radius cycle at infinity
  \(\cycle{}{\infty}=(0,0,1)\). These cycles are required for
  relations~\ref{item:quadric-flat} and~\ref{it:lobachevski-line} from
  the previous subsection. As predefined cycles, \(\cycle{}{\Space{R}{}}\) and
  \(\cycle{}{\infty}\) are placed in negative-numbered
  generations defined by the macros {{\it{}REAL\_LINE\_GEN}} and
  {{\it{}INFINITY\_GEN}}.
\item If a point is added to generation-0 of a figure, then it is
  represented by a zero-radius cycle with its centre at the given
  point. Particular parameter of such cycle dependent on the used
  metric, thus this cycle is not considered as explicitly
  defined. Thereafter, the cycle shall have some parents at a
  negative-numbered generation defined by the macro {{\it{}GHOST\_GEN}}.
\end{enumerate}
A figure can be in two different modes: {{\it{}freeze}} or {{\it{}unfreeze}},
the second is default. In the {{\it{}unfreeze}} mode an addition of a new
cycle by its relation prompts an evaluation of its parameters. If the
evaluation was successful then the obtained parameters are stored and
will be used in further calculations for all children of the cycle. Since many
relations (see the previous Subsection) are connected to quadratic
equation~\eqref{eq:det-normalisation-cond}, the solutions may come in
pairs. Furthermore, if the number or nature of conditions is not
sufficient to define the cycle uniquely (up to natural quadratic
multiplicity), then the cycle will depend on a number of free
(symbolic) variable.

There is a macro-like tool, which is called {{\bf{}subfigure}}. Such a
{{\bf{}subfigure}} is a {{\bf{}figure}} itself, such that its inner hierarchy of
generations and relations is not visible from the current
{{\bf{}figure}}. Instead, some cycles (of any generations) of the current
{{\bf{}figure}} are used as predefined cycles of generation-0 of
{{\bf{}subfigure}}. Then only one dependent cycle of {{\bf{}subfigure}}, which
is known as result, is returned back to the current {{\bf{}figure}}. The
generation of the result is calculated from generations of input
cycles by the same formula~\eqref{eq:generation-calculation}.

There is a possibility to test certain conditions (``are two cycles
orthogonal?'') or measure certain quantities (``what is their
intersection angle?'') for already defined cycles. In particular, such
methods can be used to prove geometrical statements according to the
Cartesian programme, that is replacing the synthetic geometry by
purely algebraic manipulations.
\begin{example}
  \label{ex:touch-centres-collinear}
  As an elementary demonstration, let us prove that if a cycle {{\it{}r}}
  is orthogonal to a circle {{\it{}a}} at the point {{\it{}C}} of its contact with a
  tangent line {{\it{}l}}, then {{\it{}r}} is also orthogonal to the line
  {{\it{}l}}. To simplify setup we assume that {{\it{}a}} is the unit
  circle. Here is the {\Python} code:{\upshape
  \lstset{language=Python}
\begin{lstlisting}
F=figure()
a=F.add_cycle(cycle2D(1,[0,0],-1),"a")
l=symbol("l")
C=symbol("C")
F.add_cycle_rel([is_tangent_i(a),\
  is_orthogonal(F.get_infinity()),only_reals(l)],l)
F.add_cycle_rel([is_orthogonal(C),is_orthogonal(a),\
  is_orthogonal(l),only_reals(C)],C)
r=F.add_cycle_rel([is_orthogonal(C),\
  is_orthogonal(a)],"r")
Res=F.check_rel(l,r,"cycle_orthogonal")
for i in range(len(Res)):
    print "Tangent and radius are orthogonal: %s" %\
    bool(Res[i].subs(pow(cos(wild(0)),2)\
    ==1-pow(sin(wild(0)),2)).normal())
\end{lstlisting}}
The first line creates an empty figure {{\it{}F}} with the default
euclidean metric. The next line explicitly uses parameters
\((1,0,0,-1)\) of {{\it{}a}} to add it to {{\it{}F}}. Lines~3--4 define symbols
{{\it{}l}} and {{\it{}C}}, which are needed because cycles with these labels are
defined in lines~5--8 through some relations to themselves and the
cycle {{\it{}a}}. In both cases we want to have cycles with real
coefficients only and {{\it{}C}} is additionally self-orthogonal (i.e. is a
zero-radius). Also, {{\it{}l}} is orthogonal to infinity (i.e. is a line)
and {{\it{}C}} is orthogonal to {{\it{}a}} and {{\it{}l}} (i.e. is their common
point). The tangency condition for {{\it{}l}} and self-orthogonality of
{{\it{}C}} are both quadratic relations.  The former has two solutions each
depending on one real parameter, thus line {{\it{}l}} has two
instances. Correspondingly, the point of contact {{\it{}C}} and the
orthogonal cycle {{\it{}r}} through {{\it{}C}} (defined in line~7) each have two
instances as well. Finally, lines~11--15 verify that every instance of
{{\it{}l}} is orthogonal to the respective instance of {{\it{}r}}, this is
assisted by the trigonometric substitution \(\cos^2(*)=1-\sin^2(*)\)
used for parameters of {{\it{}l}} in line~15.  The output predictably is:
\begin{verbatim}
Tangent and circle r are orthogonal: True
Tangent and circle r are orthogonal: True
\end{verbatim}
\end{example}
An original statement proved by the library {{\bf{}figure}} for the first
time will be considered in the next Section.

\section{Mathematical Usage of Libraries {{\bf{}cycle}} and {{\bf{}figure}}}
\label{sec:mathematical-results}

The developed libraries {{\bf{}cycle}} and {{\bf{}figure}} has several different uses:
\begin{itemize}
\item It is easy to produce high-quality illustrations, which are
  fully-accurate in mathematical sence. The user is not responsible
  for evaluation of cycles\vtick\ parameters, all computations are
  done by the library as soon as the figure is defined in terms of few
  geometrical relations. This is especially helpful for complicated
  images which may contain thousands of interrelated cycles. See
  Escher-like Fig.~\ref{fig:action-modular-group} which shows images
  of two circles under the modular group
  action~\cite{StewartTall02a}*{\S~14.4}.
\begin{figure}[htbp]
  \centering
  \includegraphics[,width=.9\textwidth]{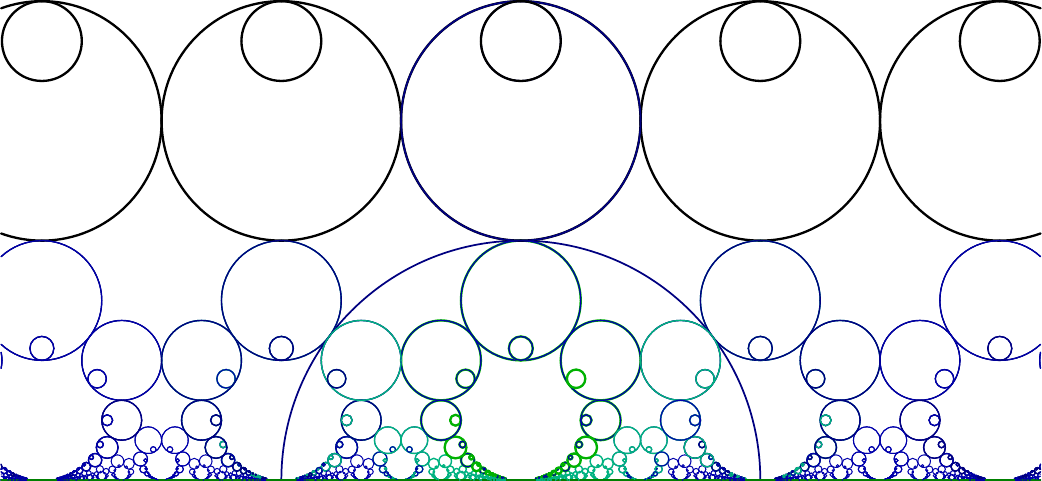}
  \caption{Action of the modular group on the upper half-plane.}
  \label{fig:action-modular-group}
\end{figure}
\begin{figure}[htbp]
  \centering
  \includegraphics[scale=.5]{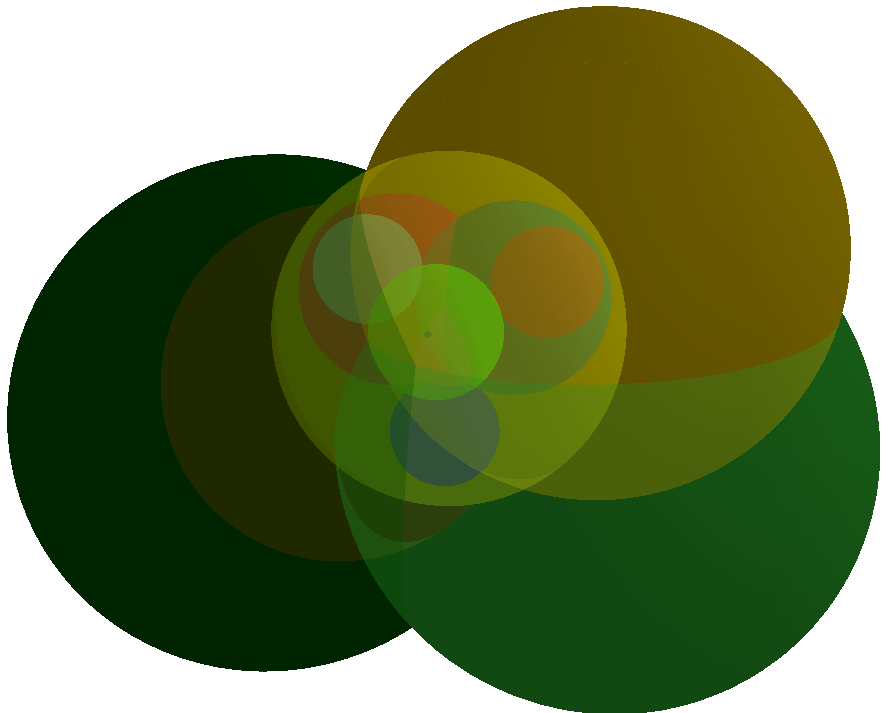}
  \caption{An example of Apollonius problem in three dimensions.}
  \label{fig:apollonius-3D}
\end{figure}
\item The package can be used for computer experiments in M\"obius--Lie
  geometry. There is a possibility to create an arrangement of cycles
  depending on one or several parameters. Then, for particular values
  of those parameters certain conditions, e.g. concurrency of cycles,
  may be numerically tested or graphically visualised. It is possible
  to create animations with gradual change of the parameters, which
  are especially convenient for illustrations, see
  Fig.~\ref{fig:nine-points-anim} and~\cite{Kisil16a}.
\item Since the library is based on the \GiNaC\ system, which provides a
  symbolic computation engine, there is a possibility to make fully
  automatic proofs of various statements in M\"obius--Lie geometry.  Usage
  of computer-supported proofs in geometry is already an established
  practice~\cites{Kisil12a,Pech07a} and it is naturally to expect its further
  rapid growth.
\item Last but not least, the combination of classical beauty of Lie
  sphere geometry and modern computer technologies is a useful
  pedagogical tool to widen interest in mathematics through visual and
  hands-on experience.
\end{itemize}

Computer experiments are especially valuable for Lie geometry of
indefinite or nilpotent metrics since our intuition is not elaborated
there in contrast to the Euclidean
space~\cites{Kisil07a,Kisil06a,Kisil05a}. Some advances in
the two-dimensional space were achieved
recently~\cites{Mustafa17a,Kisil12a}, however further developments in
higher dimensions are still awaiting their researchers.

\begin{figure}[htbp]
  \centering
  \makebox[0pt][l]{(a)}\includegraphics[scale=.58]{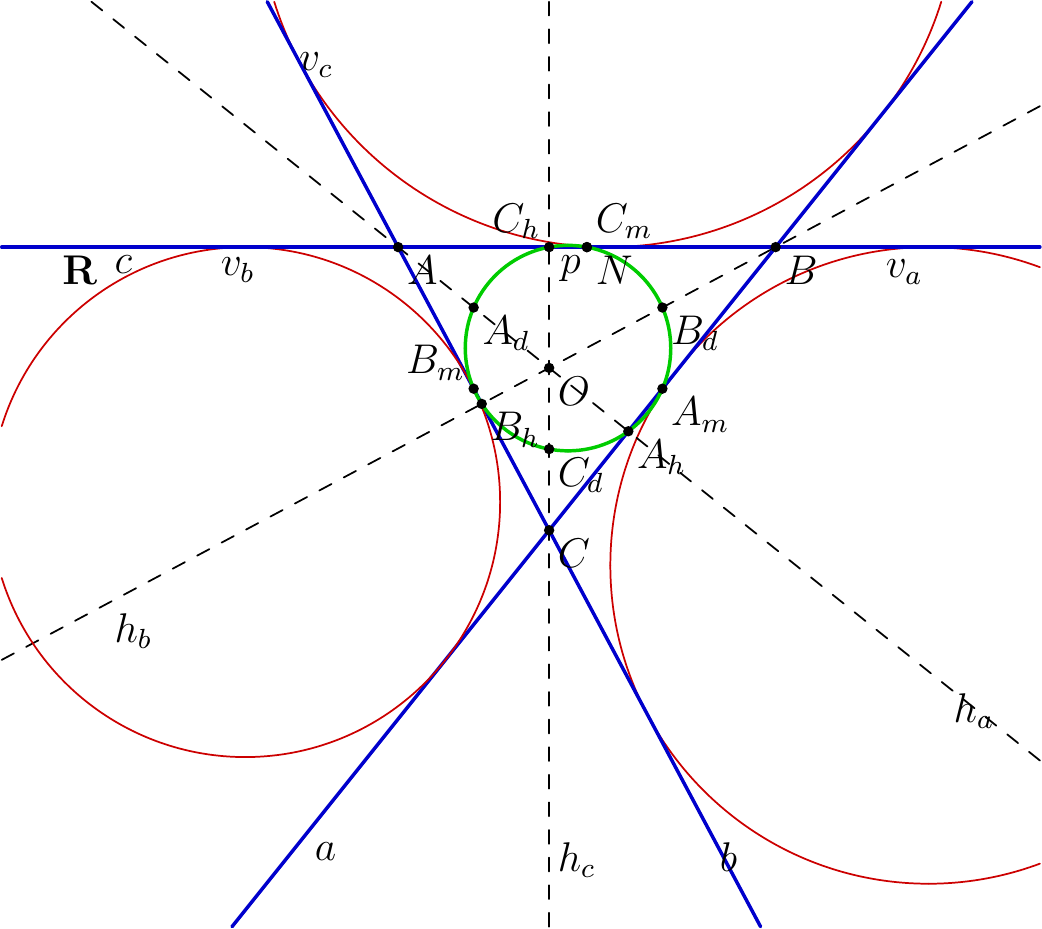}\hfill
  \makebox[0pt][l]{(b)}\includegraphics[scale=.58]{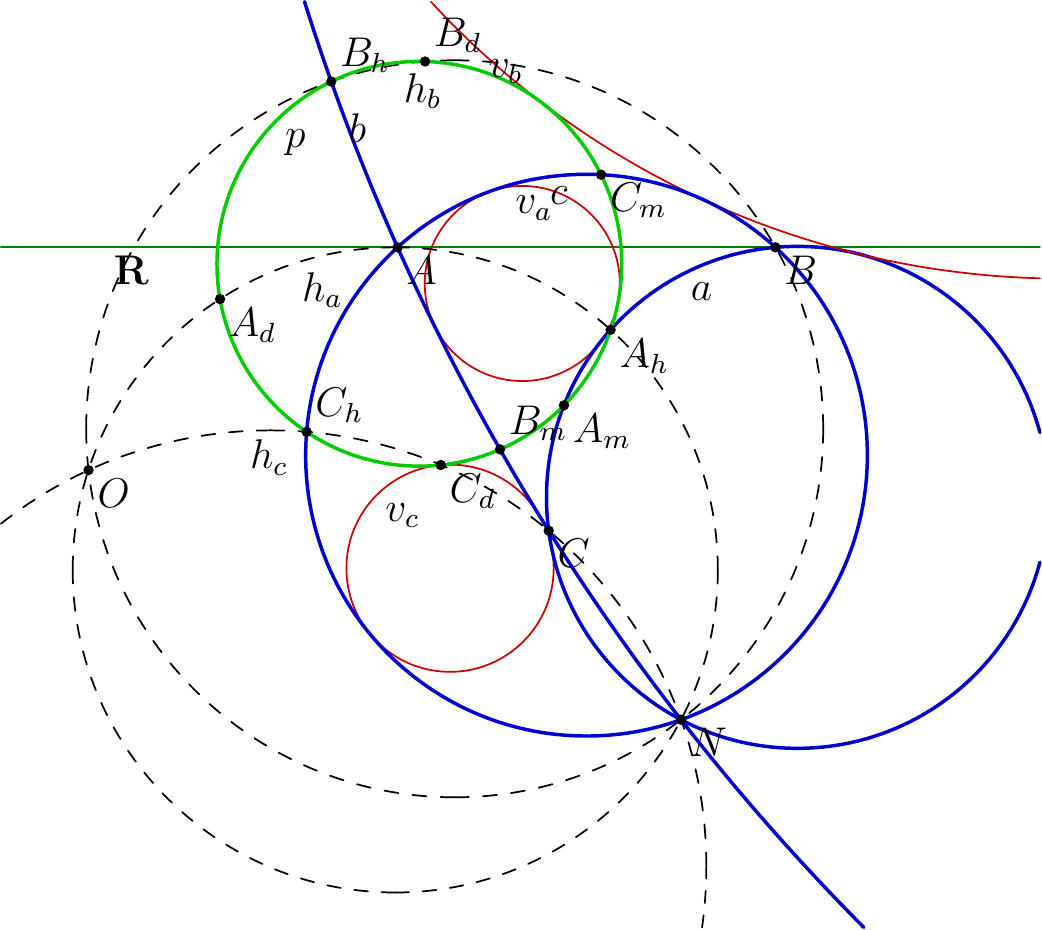}\\[1em]
  \makebox[0pt][l]{(c)}\includegraphics[scale=.58]{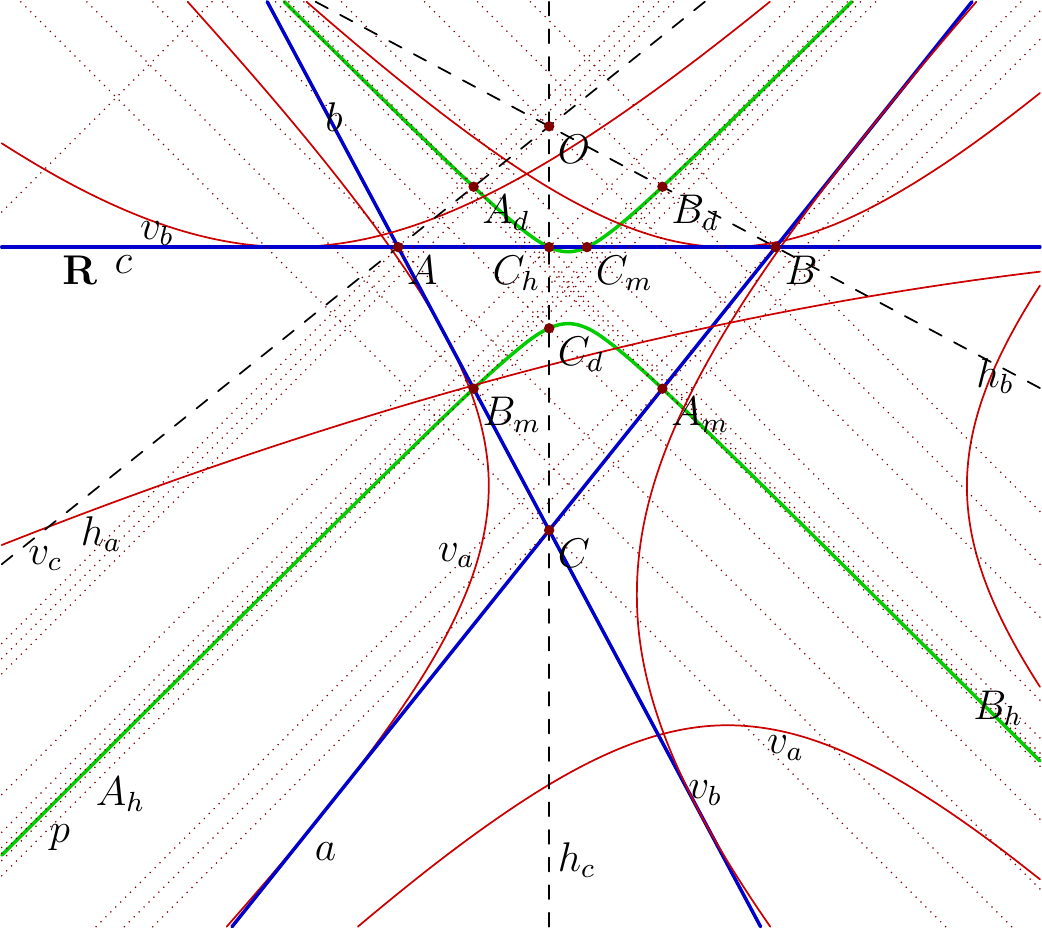}\hfill
  \makebox[0pt][l]{(d)}\includegraphics[scale=.58]{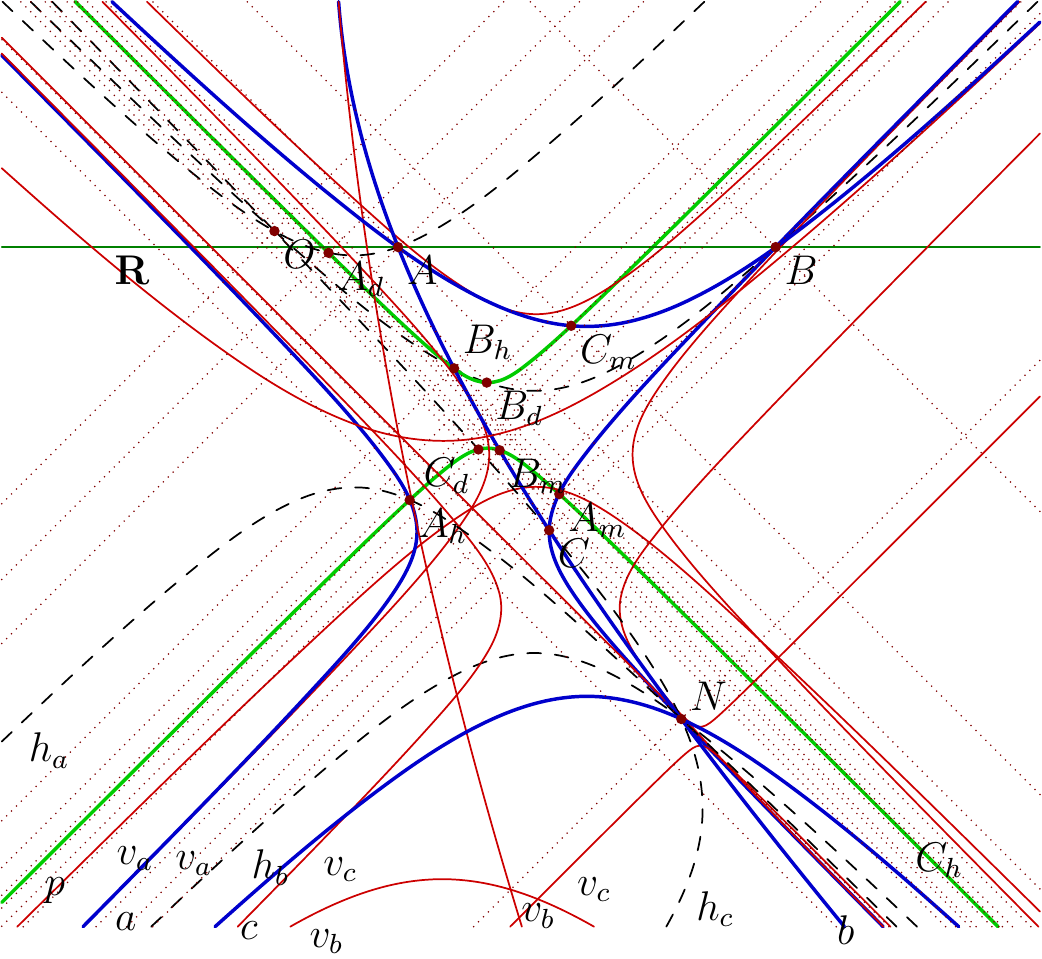}
  \caption[The illustration of the conformal nine-points theorem]
  {The illustration of the conformal nine-points theorem. The
    left column is the statement for a triangle with straight sides
    (the point {{\it{}N}} is at infinity), the right column is its
    conformal version (the point {{\it{}N}} is at the finite part). The
    first row show the elliptic point space, the second row---the
    hyperbolic point space. Thus, the top-left picture shows the
    traditional theorem, three other pictures---its different modifications.}
  \label{fig:illustr-conf-nine}
\end{figure}

As a non-trivial example of automated proof accomplished by the {{\bf{}figure}}
library for the first time, we present a FLT-invariant version of the
classical nine-point theorem~\citelist{\cite{Pedoe95a}*{\S~I.1}
  \cite{CoxeterGreitzer}*{\S~1.8}},
cf. Fig.~\ref{fig:illustr-conf-nine}(a):
\begin{theorem}[Nine-point cycle]
  \label{th:nine-points}
  Let \(ABC\) be an arbitrary triangle with the orthocenter (the
  points of intersection of three altitudes) \(H\), then
  the following nine points belongs to the same cycle, which may be a
  circle or a hyperbola:
  \begin{enumerate}
    \item Foots of three altitudes, that is points of pair-wise
      intersections \(AB\) and \(CH\), \(AC\) and \(BH\), \(BC\) and \(AH\).
    \item Midpoints of sides \(AB\), \(BC\) and \(CA\).
    \item Midpoints of intervals \(AH\), \(BH\) and \(CH\).
  \end{enumerate}
\end{theorem}
There are many further interesting properties,
e.g. nine-point circle is externally tangent to that triangle three
excircles and internally tangent to its incircle as it seen from
Fig.~\ref{fig:illustr-conf-nine}(a).

To adopt the statement for cycles geometry we need to find a
FLT-invariant meaning of the midpoint \(A_m\) of an interval \(BC\),
because the equality of distances \(BA_m\) and \(A_mC\) is not
FLT-invariant. The definition in cycles geometry can be done by either
of the following equivalent relations:
\begin{itemize}
\item The midpoint \(A_m\) of an interval \(BC\) is defined by the
  cross-ratio \(\frac{BA_m}{CA_m} : \frac{BI}{CI}=1\), where \(I\) is
  the point at infinity.
\item We construct the midpoint \(A_m\) of an interval \(BC\) as the
  intersection of the interval and the line orthogonal to \(BC\) and
  to the cycle, which uses \(BC\) as its diameter. The latter
  condition means that the cycle passes both points \(B\) and \(C\)
  and is orthogonal to the line \(BC\).
\end{itemize}
Both procedures are meaningful if we replace the point at infinity
\(I\) by an arbitrary fixed point \(N\) of the plane. In the second
case all lines will be replaced by cycles passing through \(N\), for
example the line through \(B\) and \(C\) shall be replaced by a cycle
through \(B\), \(C\) and \(N\). If we similarly replace ``lines'' by
``cycles passing through \(N\)'' in Thm.~\ref{th:nine-points} it turns
into a valid FLT-invariant version,
cf. Fig.~\ref{fig:illustr-conf-nine}(b). Some additional properties,
e.g. the tangency of the nine-points circle to the ex-/in-circles, are
preserved in the new version as well.  Furthermore, we can illustrate
the connection between two versions of the theorem by an animation,
where the infinity is transformed to a finite point \(N\) by a
continuous one-parameter group of FLT,
see. Fig.~\ref{fig:nine-points-anim} and further examples
at~\cite{Kisil16a}.

It is natural to test the nine-point theorem in the hyperbolic and the
parabolic spaces. Fortunately, it is very easy under the given
implementation: we only need to change the defining metric of the
point space, cf.~\cite{Allen41a}, this can be done for an already
defined figure. The corresponding figures
Fig.~\ref{fig:illustr-conf-nine}(c) and~(d) suggest that the
hyperbolic version of the theorem is still true in the plain and even
FLT-invariant forms. We shall clarify that the hyperbolic version of
the Thm.~\ref{th:nine-points} specialises the nine-point conic of a
complete quadrilateral \cites{CerinGianella06a,DeVilliers06a}: in
addition to the existence of this conic, our theorem specifies its
type for this particular arrangement as equilateral hyperbola with the
vertical axis of symmetry.

\begin{figure}[htbp]
  \centering
  \animategraphics[controls=true,width=.9\textwidth]{50}{_nine-points-anim}{}{}
  \caption[Animated transition between the classical and conformal
  nine-point theorems]{Animated transition between the classical and
    conformal versions of the nine-point theorem. Use control buttons to activate
    it. You may need \textsf{Adobe Acrobat Reader} for this feature.}
  \label{fig:nine-points-anim}
\end{figure}

The computational power of the package is sufficient not only to hint
that the new theorem is true but also to make a complete proof. To
this end we define an ensemble of cycles with exactly same
interrelations, but populate the generation-0 with points \(A\), \(B\)
and \(C\) with symbolic coordinates, that is, objects of the \GiNaC\
{{\bf{}class}\ {\bf{}realsymbol}}. Thus, the entire figure defined from them will
be completely general. Then, we may define the hyperbola passing
through three bases of altitudes and check by the symbolic
computations that this hyperbola passes another six ``midpoints'' as
well.

In the parabolic space the nine-point Thm.~\ref{th:nine-points} is not
preserved in this manner. It is already
observed~\cites{Kisil12a,Kisil05a,%
  Kisil15a,Kisil07a,Kisil09e,Kisil11a,Mustafa17a,BarrettBolt10a}, that
the degeneracy of parabolic metric in the point space requires certain
revision of traditional definitions. The parabolic variation of
nine-point theorem may prompt some further considerations as well. An
expanded discussion of various aspects of the nine-point construction
shall be the subject of a separate paper.



\section*{Acknowledgement}
\label{sec:acknowledgement}

I am grateful to Prof.~Jay~P.~Fillmore for stimulating discussion,
which enriched the library {{\bf{}figure}}.
The University of Leeds provided a generous summer internship to work
on Graphical User Interface to the library, which was initiated by
Luke Hutton with skills and enthusiasm. Cameron Kumar wrote the
initial version of a
\href{https://sourceforge.net/projects/cycle3dvis.moebinv.p/}{3D cycle
  visualiser} as a part of his BSc project at the University of Leeds.

\providecommand{\noopsort}[1]{} \providecommand{\printfirst}[2]{#1}
  \providecommand{\singleletter}[1]{#1} \providecommand{\switchargs}[2]{#2#1}
  \providecommand{\irm}{\textup{I}} \providecommand{\iirm}{\textup{II}}
  \providecommand{\vrm}{\textup{V}} \providecommand{\cprime}{'}
  \providecommand{\eprint}[2]{\texttt{#2}}
  \providecommand{\myeprint}[2]{\texttt{#2}}
  \providecommand{\arXiv}[1]{\myeprint{http://arXiv.org/abs/#1}{arXiv:#1}}
  \providecommand{\doi}[1]{\href{http://dx.doi.org/#1}{doi:
  #1}}\providecommand{\CPP}{\texttt{C++}}
  \providecommand{\NoWEB}{\texttt{noweb}}
  \providecommand{\MetaPost}{\texttt{Meta}\-\texttt{Post}}
  \providecommand{\GiNaC}{\textsf{GiNaC}}
  \providecommand{\pyGiNaC}{\textsf{pyGiNaC}}
  \providecommand{\Asymptote}{\texttt{Asymptote}}

\end{document}